\documentclass[11pt,draft]{amsart}
\usepackage{amscd}
\usepackage{amssymb}

\newtheorem{thm}[equation]{Theorem}
\newtheorem{cor}[equation]{Corollary}
\newtheorem{prop}[equation]{Proposition}
\newtheorem{lem}[equation]{Lemma}
\theoremstyle{definition}
\newtheorem{dfn}[equation]{Definition}
\newtheorem{rem}[equation]{Remark}
\newtheorem{exa}[equation]{Example}

\newtheorem{prob}[equation]{Problem}

\numberwithin{equation}{section}

\newcommand{\iso}{\stackrel{\simeq}{\rightarrow}}
\newcommand{\inj}{\hookrightarrow}
\newcommand{\surj}{\twoheadrightarrow}
\newcommand{\ar}{\rightarrow}
\newcommand{\lar}{\leftarrow}
\newcommand{\opn}{\operatorname}
\newcommand{\bdot}{{\textstyle \cdot}}
\newcommand{\ul}{\underline}
\newcommand{\blnk}[1]{\mbox{\hspace{#1}}}
\newcommand{\rmitem}[1]{\item[\text{\textup{(#1)}}]}
\newcommand{\mfrak}[1]{\mathfrak{#1}}
\newcommand{\mcal}[1]{\mathcal{#1}}
\newcommand{\msf}[1]{\mathsf{#1}}
\newcommand{\mbf}[1]{\mathbf{#1}}
\newcommand{\mrm}[1]{\mathrm{#1}}
\newcommand{\mbb}[1]{\mathbb{#1}}
\newcommand{\gfrac}[2]{\genfrac{[}{]}{0pt}{}{#1}{#2}}
\title[Residues on Schemes]{Residues and Differential
Operators on Schemes}
\author{Amnon Yekutieli}
\address{Department of Theoretical Mathematics,
The Weizmann Institute of Science,
Rehovot 76100, ISRAEL}
\date{3.6.97}
\email{amnon@wisdom.weizmann.ac.il}
\thanks{The author is supported by an Allon Fellowship, and is the
incumbent of the Anna and Maurice Boukstein Career Development Chair.}

\begin{document}
\maketitle
\tableofcontents

\section{Introduction}

Suppose $X$ is a finite type scheme over a field $k$, with structural
morphism $\pi$. Consider the twisted inverse image functor
$\pi^{!} : \msf{D}^{+}_{\mrm{c}}(k) \ar
\msf{D}^{+}_{\mrm{c}}(X)$
of Grothendieck Duality Theory (see \cite{RD}).
The {\em residue complex} $\mcal{K}^{\bdot}_{X}$
is defined to be the Cousin complex of $\pi^{!} k$.
It is a bounded complex of quasi-coherent $\mcal{O}_{X}$-modules,
possessing remarkable functorial properties.
In this paper we provide an explicit construction of
$\mcal{K}^{\bdot}_{X}$. This construction reveals some new properties
of $\mcal{K}^{\bdot}_{X}$, and also has applications in other
areas of algebraic geometry.

Grothendieck Duality, as developed by Hartshorne in \cite{RD},
is an abstract theory, stated in the language of derived categories.
Even though this abstraction is suitable for many important
applications, often one wants more explicit information.
Thus a significant
amount of work was directed at finding a presentation of
duality in terms of differential forms and residues.
Mostly the focus was on the dualizing sheaf $\omega_{X}$,
in various circumstances. The structure of $\omega_{X}$ as a
coherent $\mcal{O}_{X}$-module and its variance properties are
thoroughly understood by now, thanks to an extended effort
including \cite{Kl}, \cite{KW}, \cite{Li}, \cite{HK1}, \cite{HK2},
\cite{LS} and \cite{HS}.
Regarding an explicit presentation of the full duality theory of
dualizing complexes, there have been some advances in recent years,
notably in the papers \cite{Ye1}, \cite{SY}, \cite{Hu}, \cite{Hg}
and \cite{Sa}.

In this paper we give a totally new construction of the residue
complex $\mcal{K}^{\bdot}_{X}$, when $k$ is a perfect field of
any characteristic and $X$ is any finite type $k$-scheme.
The main idea is the use of {\em Beilinson Completion Algebras}
(BCAs), which were introduced in \cite{Ye2}.
These algebras are generalizations of complete local rings, and they
carry a mixed algebraic-analytic structure.
A review of BCAs and their properties is included in Section 1,
for the reader's convenience.

Given a point $x \in X$, the complete local ring
$\widehat{\mcal{O}}_{X, x} = \mcal{O}_{X, (x)}$
is a BCA, so according to \cite{Ye2} it has a {\em dual module}
$\mcal{K}(\mcal{O}_{X, (x)})$. This module is a canonical model
for the injective hull of the residue field $k(x)$.
If $(x,y)$ is a saturated chain of points (i.e.\ $y$ is an immediate
specialization of $x$) then there is a BCA
$\mcal{O}_{X, (x,y)}$ and homomorphisms
$\mrm{q} : \mcal{K}(\mcal{O}_{X, (x)}) \ar
\mcal{K}(\mcal{O}_{X, (x,y)})$
and
$\opn{Tr} : \mcal{K}(\mcal{O}_{X, (x,y)}) \ar
\mcal{K}(\mcal{O}_{X, (y)})$.
The dual modules $\mcal{K}(-)$ and the homomorphisms
$\mrm{q}$ and $\opn{Tr}$ have explicit formulas in terms of differential
forms and coefficient fields. Set
$\delta_{(x,y)} := \opn{Tr} \mrm{q} : \mcal{K}(\mcal{O}_{X, (x)}) \ar
\mcal{K}(\mcal{O}_{X, (y)})$.
Define a graded quasi-coherent sheaf $\mcal{K}_{X}^{\bdot}$ by
\[ \mcal{K}_{X}^{q} := \bigoplus_{\opn{dim} \overline{\{x\}} = -q}
\mcal{K}(\mcal{O}_{X, (x)}) \]
and a degree $1$ homomorphism
\[ \delta := (-1)^{q+1} \sum_{(x,y)} \delta_{(x,y)} . \]

It turns out that
$(\mcal{K}_{X}^{\bdot}, \delta)$ is a residual complex on $X$, and
it is canonically isomorphic to $\pi^{!} k$ in the derived category
$\msf{D}(X)$. Hence it is the residue complex of $X$, as defined in
the first paragraph. The functorial
properties of $\mcal{K}_{X}^{\bdot}$ w.r.t.\ proper and \'{e}tale
morphisms are obtained directly from
corresponding properties of BCAs, and therefore are reduced to
explicit formulas. All this is worked out in Sections 2 and 3.

An $\mcal{O}_{X}$-module $\mcal{M}$ has a dual complex
$\opn{Dual} \mcal{M} := \mcal{H}om^{\bdot}_{\mcal{O}_{X}}(\mcal{M},
\mcal{K}_{X}^{\bdot})$.
Suppose
$\mrm{d} : \mcal{M} \ar \mcal{N}$
is a differential operator (DO). In Theorem \ref{thm3.1} we prove
there is a dual operator
$\opn{Dual}(\mrm{d}) : \opn{Dual} \mcal{N} \ar \opn{Dual} \mcal{M}$,
which commutes with $\delta$. The existence of $\opn{Dual}(\mrm{d})$
does not follow from
formal considerations of duality theory; it is a consequence of our
particular construction using BCAs
(but cf.\ Remarks \ref{rem3.2} and \ref{rem3.3}).
The construction also provides explicit formulas for
$\opn{Dual}(\mrm{d})$ in terms of differential operators and
residues, which are used in the applications in Sections 6 and 7.

Suppose $A$ is a finite type $k$-algebra, and let
$\mcal{D}(A)$ be the ring of differential operators of $A$.
As an immediate application of Theorem \ref{thm3.1} we obtain a
description of the opposite ring $\mcal{D}(A)^{\circ}$,
as the ring of DOs on
$\mcal{K}_{A}^{\bdot}$ which commute with $\delta$
(Theorem \ref{thm3.2}). In the case of a
Gorenstein algebra it follows that the opposite ring
$\mcal{D}(A)^{\circ}$ is naturally isomorphic to
$\omega_{A} \otimes_{A} \mcal{D}(A) \otimes_{A} \omega_{A}^{-1}$
(Corollary \ref{cor3.6}).

Applying Theorem \ref{thm3.1} to the De Rham complex
$\Omega^{\bdot}_{X/k}$ we obtain the {\em De Rham-residue complex}
$\mcal{F}_{X}^{\bdot} = \opn{Dual} \Omega^{\bdot}_{X/k}$.
Up to signs this coincides with El-Zein's complex
$\mcal{K}_{X}^{\bdot,*}$ of \cite{EZ} (Corollary \ref{cor4.3}).
The fundamental class $\mrm{C}_{Z} \in \mcal{F}_{X}^{\bdot}$,
for a closed subscheme $Z \subset X$, is easily described in this context
(Definition \ref{dfn4.4}).

The construction above works also for a formal scheme
$\mfrak{X}$ which is of formally finite type over $k$, in the sense
of \cite{Ye3}. An example of such a formal scheme is the completion
$\mfrak{X} = Y_{/X}$, where $X$
is a locally closed subset of the finite type $k$-scheme $Y$.
Therefore we get a complex
$\mcal{F}_{\mfrak{X}}^{\bdot} =
\opn{Dual} \widehat{\Omega}^{\bdot}_{\mfrak{X}/k}$.
When $X \subset \mfrak{X}$ is a smooth formal embedding
(see Definition \ref{dfn5.3}) we prove that the cohomology modules
$\mrm{H}^{q}(X, \mcal{F}_{\mfrak{X}}^{\bdot})$
are independent of $\mfrak{X}$. This is done by analyzing the
$E_{1}$ term of the {\em niveau spectral sequence} converging to
$\mrm{H}^{\bdot}(X, \mcal{F}_{\mfrak{X}}^{\bdot})$
(Theorem \ref{thm5.3}). Here we assume $\opn{char} k = 0$.
The upshot is that
$\mrm{H}^{q}(X, \mcal{F}_{\mfrak{X}}^{\bdot}) =
\mrm{H}_{-q}^{\mrm{DR}}(X)$,
the De Rham homology.
There is an advantage in using smooth formal embeddings.
If $U \ar X$ is any \'{e}tale morphism, then there
is an \'{e}tale morphism $\mfrak{U} \ar \mfrak{X}$ s.t.\
$U = \mfrak{U} \times_{\mfrak{X}} X$, so
$U \subset \mfrak{U}$ is a smooth formal embedding. From this we
conclude that
$\mrm{H}_{\bdot}^{\mrm{DR}}(-)$
is a contravariant functor on $X_{\mrm{et}}$, the small \'{e}tale
site. Previously it was only known that
$\mrm{H}_{\bdot}^{\mrm{DR}}(-)$ is contravariant for open
immersions (cf.\ \cite{BlO} Example 2.2).

Suppose $X$ is smooth, and let
$\mcal{H}^{p}_{\mrm{DR}}$
be the sheafification of the presheaf
$U \mapsto \mrm{H}^{p}_{\mrm{DR}}(U)$ on $X_{\mrm{Zar}}$.
Bloch-Ogus \cite{BlO} give a flasque resolution of
$\mcal{H}^{p}_{\mrm{DR}}$, the {\em arithmetic resolution}.
It involves the sheaves
$i_{x*} \mrm{H}^{q} \Omega^{\bdot}_{k(x)/k}$
where $i_{x} : \{x\} \ar X$ is the inclusion map. Our analysis of the
niveau spectral sequence shows that the coboundary operator of this
resolution is a sum of Parshin residues (Corollary \ref{cor5.2}).

Our final application of the new construction of the residue complex
is to describe the {\em intersection cohomology $\mcal{D}$-module}
$\mcal{L}(X, Y)$, when $X$ is an integral curve embedded in a smooth
$n$-dimensional variety $Y$ (see \cite{BrKa}). Again we assume $k$
has characteristic $0$.
In fact we are able to describe all coherent $\mcal{D}_{Y}$-submodules
of $\mcal{H}^{n-1}_{X} \mcal{O}_{Y}$
in terms of the singularities of $X$ (Corollary \ref{cor6.10}).
This description is an algebraic version of Vilonen's work in \cite{Vi},
replacing complex analysis with BCAs and algebraic residues.
It is our hope that a similar description will be found in the general
case, namely $\opn{dim} X > 1$. Furthermore, we hope to give in the
future an explicit description of the Cousin complex of
$\opn{DR} \mcal{L}(X, Y) = \Omega^{\bdot}_{Y / k} \otimes
\mcal{L}(X, Y)$.
Note that for $X=Y$ one has
$\mcal{L}(X, Y) = \mcal{O}_{X}$, so this Cousin complex is
nothing but $\mcal{F}^{\bdot}_{X}$.

\medskip \noindent
{\bf Acknowledgments.}\
I wish to thank J.\ Lipman and S.\ Kleiman for their continued interest
in this work. Thanks also to
P.\ Sastry, R.\ H\"{u}bl, V.\ Hinich, V.\ Berkovich, H.\ Esnault,
K.\ Smith and K.\ Vilonen for helpful conversations, and
thanks to the referee for valuable advice on Sections 4 and 7.

\section{Review of Beilinson Completion Algebras}

Let us begin by reviewing some facts about Topological Local Fields
(TLFs) and Beilinson Completion Algebras (BCAs) from the papers
\cite{Ye1} and \cite{Ye2}.

A semi-topological (ST) ring is a ring $A$, with a linear topology on
its underlying additive group, such that for every $a \in A$ the
multiplication (on either side)
$a : A \ar A$ is continuous.

Let $K$ be a field. We say $K$ is an $n$-dimensional local field
if there is a sequence of complete discrete valuation rings
$\mcal{O}_{1}, \ldots, \mcal{O}_{n}$,
where the fraction field of $\mcal{O}_{1}$ is $K$,
and the residue field of $\mcal{O}_{i}$ is
the fraction field of $\mcal{O}_{i + 1}$.

Fix a perfect field $k$.
A {\em topological local field} of dimension $n$ over $k$ is a
$k$-algebra
$K$ with structures of semi-topological ring and $n$-dimensional local
field, satisfying the following parameterization condition:
there exists an isomorphism of $k$-algebras
$K \cong F((s_{1}, \ldots, s_{n}))$
for some field $F$, finitely generated over $k$,
which respects the two structures.
Here
$F((s_{1}, \ldots, s_{n})) = F((s_{n})) \cdots ((s_{1}))$
is the field of iterated Laurent series, with its inherent topology
and valuation rings ($F$ is discrete).
One should remark that for $n = 1$ we are in the classical situation,
whereas for $n \geq 2$, $F((s_{1}, \ldots, s_{n}))$
is not a topological ring.

TLFs make up a category $\msf{TLF}(k)$,
where a morphism $K \ar L$ is a continuous
$k$-algebra homomorphism which preserves the valuations, and the induced
homomorphism of the last residue fields is finite.
Write $\Omega^{\mrm{sep}, \bdot}_{K / k}$ for the separated
algebra of differentials; with the parameterization above
$\Omega^{\mrm{sep}, \bdot}_{K / k} \cong
K \otimes_{F} \Omega^{\bdot}_{F[\, \ul{s}\, ] / k}$.
Then there is a functorial residue map
$\opn{Res}_{L / K} : \Omega^{\mrm{sep}, \bdot}_{L / k} \ar
\Omega^{\mrm{sep}, \bdot}_{K / k}$
which is $\Omega^{\mrm{sep}, \bdot}_{K / k}$-linear and
lowers degree by $\opn{dim} L / K$. For instance if
$L = K((t))$ then
\begin{equation} \label{eqn10.1}
\opn{Res}_{L / K} \left( \sum_{i} t^{i} \mrm{d} t \wedge \alpha_{i}
\right) =
\alpha_{-1} \in \Omega^{\mrm{sep}, \bdot}_{K / k} .
\end{equation}
TLFs and residues were initially developed by Parshin and Lomadze,
and the theory was enhanced in \cite{Ye1}.

The notion of a {\em Beilinson completion algebra}
was introduced in \cite{Ye2}. A BCA is a semi-local,
semi-topological $k$-algebra, each of whose residue fields
$A / \mfrak{m}$ is a topological local field. Again there is a
parameterization condition: when $A$ is local,
there should exist a surjection
\[ F((\ul{s})) [[\, \ul{t}\, ]] =
F((\ul{s})) [[ t_{1}, \ldots, t_{m} ]] \surj A \]
which is strict topologically (i.e.\ $A$ has the quotient topology)
and respects the valuations.
Here $F((\ul{s}))$ is as above and
$F((\ul{s})) [[\, \ul{t}\, ]]$ is the ring of formal power series over
$F((\ul{s}))$.
The notion of BCA is an abstraction of the algebra gotten by
Beilinson's completion, cf.\ Lemma \ref{lem10.1}.

There are two distinguished kinds of homomorphisms between BCAs.
The first kind is a {\em morphism of BCAs} $f : A \ar B$
(see \cite{Ye2} Definition 2.5),
and the category they constitute is denoted $\msf{BCA}(k)$.
A morphism is continuous, respects the valuations on the residue fields,
but in general is not a local homomorphism.
For instance, the homomorphisms
$k \ar k[[ s, t ]] \ar k((s))[[t]] \ar k((s))((t))$ are all morphisms.
$\msf{TLF}(k)$ is a full subcategory of $\msf{BCA}(k)$, consisting
of those BCAs which are fields.

The second kind of homomorphism is an {\em intensification homomorphism}
$u : A \ar \widehat{A}$ (see \cite{Ye2} Definition 3.6).
An intensification is flat, topologically \'{e}tale
(relative to $k$) and unramified (in the appropriate sense).
It can be viewed as a sort of localization or
completion. Here examples are
$k(s)[[t]] \ar k((s))[[t]]$ and
$k(s,t) \ar k(s)((t)) \ar k((s))((t))$.

Suppose $f : A \ar B$ is a morphism of BCAs and
$u : A \ar \widehat{A}$ is an intensification. The Intensification Base
Change Theorem (\cite{Ye2} Theorem 3.8) says there is a BCA
$\widehat{B} = B \otimes_{A}^{(\wedge)} \widehat{A}$, a morphism
$\widehat{f} : \widehat{A} \ar \widehat{B}$ and an intensification
$v : B \ar \widehat{B}$, with
$v f = \widehat{f} u$. These are determined up to isomorphism and
satisfy certain universal properties.
For instance,
$k((s))[[t]] = k(s)[[t]] \otimes_{k(s)}^{(\wedge)} k((s))$.

According to  \cite{Ye2} Theorem 6.14, every $A \in \msf{BCA}(k)$
has a dual module $\mcal{K}(A)$.
The module $\mcal{K}(A)$ is a ST $A$-module. Algebraically it is
an injective hull of $A / \mfrak{r}$, where $\mfrak{r}$ is
the Jacobson radical.
$\mcal{K}(A)$ is also a right $\mcal{D}(A)$-module, where
$\mcal{D}(A)$ denotes the ring of continuous differential operators
of $A$ (relative to $k$).
For a ST $A$-module $M$ let
$\opn{Dual}_{A} M := \opn{Hom}_{A}^{\mrm{cont}}(M, \mcal{K}(A))$.

The dual modules have variance properties w.r.t.\ morphisms and
intensifications.
Given a morphism of BCAs $f : A \ar B$, according to \cite{Ye2}
Theorem 7.4 there is an $A$-linear map
$\opn{Tr}_{f} : \mcal{K}(B) \ar \mcal{K}(A)$.
This induces an isomorphism
$\mcal{K}(B) \cong
\opn{Hom}_{A}^{\mrm{cont}}(B, \mcal{K}(A))$.
Given an intensification homomorphism $u : A \ar \widehat{A}$,
according to \cite{Ye2} Proposition 7.2 there is an $A$-linear map
$\mrm{q}_{u} : \mcal{K}(A) \ar \mcal{K}(\widehat{A})$.
It induces an isomorphism
$\mcal{K}(\widehat{A}) \cong \widehat{A} \otimes_{A} \mcal{K}(A)$.
Furthermore $\opn{Tr}$ and $\mrm{q}$ commute across intensification
base change:
$\opn{Tr}_{\widehat{B} / \widehat{A}} \mrm{q}_{B / \widehat{B}} =
\mrm{q}_{\widehat{A} / A} \opn{Tr}_{B / A}$.

In case of a TLF $K$, one has
$\mcal{K}(K) = \omega(K) = \Omega^{p, \mrm{sep}}_{K/k}$,
where $p = \opn{rank} \Omega^{1, \mrm{sep}}_{K/k}$.
For a morphism of TLFs $f :K \ar L$ one has
$\opn{Tr}_{f} = \opn{Res}_{f}$, whereas for an intensification
$u : K \ar \widehat{K}$ the homomorphism
$\mrm{q}_{u} : \Omega^{p, \mrm{sep}}_{K/k} \ar
\Omega^{p, \mrm{sep}}_{\widehat{K}/k}$
is the canonical inclusion for a topologically \'{e}tale extension
of fields.

\begin{exa} \label{exa10.1}
Take
$L := k(s,t)$, $\widehat{L} := k(s)((t))$,
$A := k(s)[[t]]$,
$\widehat{A} := k((s))[[t]]$,
$K := k(s)$
and
$\widehat{K} := k((s))$.
The inclusions
$L \ar \widehat{L}$, $K \ar \widehat{K}$ and
$A \ar \widehat{A}$
are intensifications, whereas
$K \ar A \ar \widehat{L}$ and $\widehat{K} \ar \widehat{A}$
are morphisms. Using the isomorphism
$\mcal{K}(A) \cong
\opn{Hom}_{K}^{\mrm{cont}}(A, \Omega^{1, \mrm{sep}}_{K / k})$
induced by $\opn{Tr}_{A / K}$,
we see that for $\alpha \in \Omega^{2, \mrm{sep}}_{\widehat{L} / k}$
the element
$\opn{Tr}_{\widehat{L} / A}(\alpha) \in \mcal{K}(A)$
is represented by the functional
$a \mapsto \opn{Res}_{\widehat{L} / K}(a \alpha)$,
$a \in A$.
Also for
$\phi \in \mcal{K}(A)$ the element
$\widehat{\phi} = \mrm{q}_{\widehat{A} / A}(\phi) \in
\mcal{K}(\widehat{A})$
is represented by the unique $\widehat{K}$-linear functional
$\widehat{\phi} : \widehat{A} \ar
\Omega^{1, \mrm{sep}}_{\widehat{K} / k}$
extending $\phi$.
\end{exa}

\begin{rem}
The proof of existence of dual modules with their variance properties
in \cite{Ye2} is straightforward,
using Taylor series expansions, differential operators
and the residue pairing.
\end{rem}

Let $A$ be a noetherian ring and $\mfrak{p}$ a prime ideal.
Consider the exact functor on $A$-modules
$M \mapsto M_{(\mfrak{p})} := \widehat{A}_{\mfrak{p}} \otimes_{A} M$.
For $M$ finitely generated we have
$M_{(\mfrak{p})} \cong
\lim_{\lar i} M_{\mfrak{p}} / \mfrak{p}^{i} M_{\mfrak{p}}$,
and if $M = \lim_{\alpha \ar} M_{\alpha}$, then
$M_{(\mfrak{p})} \cong
\lim_{\alpha \ar} (M _{\alpha})_{(\mfrak{p})}$.
This was generalized by Beilinson (cf.\ \cite{Be}) as follows.

\begin{dfn} \label{dfn10.1}
Let $M$ be an $A$-module and let
$\xi = (\mfrak{p}_{0}, \ldots, \mfrak{p}_{n})$ be a chain of prime
ideals, namely $\mfrak{p}_{i} \subset \mfrak{p}_{i + 1}$.
Define the {\em Beilinson completion} $M_{\xi}$
by recursion on $n$, $n \geq -1$.
\begin{enumerate}
\item If $n = -1$ (i.e.\ $\xi = \emptyset$),
let $M_{\xi} := M$ with the discrete topology.
\item If $n \geq 0$ and $M$ is finitely generated, let
\[ M_{(\mfrak{p}_{0}, \ldots, \mfrak{p}_{n})} :=
\lim_{\lar i}\, (M_{\mfrak{p}_{0}} / \mfrak{p}_{0}^{i}
M_{\mfrak{p}_{0}})_{(\mfrak{p}_{1}, \ldots, \mfrak{p}_{n})} . \]
\item For arbitrary $M$, let $\{ M_{\alpha} \}$ be the set of
finitely generated submodules of $M$, and let
\[ M_{(\mfrak{p}_{0}, \ldots, \mfrak{p}_{n})} :=
\lim_{\alpha \ar}\, (M_{\alpha})_{(\mfrak{p}_{0}, \ldots, \mfrak{p}_{n})}
. \]
\end{enumerate}
\end{dfn}

A chain $\xi = (\mfrak{p}_{0}, \ldots, \mfrak{p}_{n})$
is {\em saturated} if $\mfrak{p}_{i + 1}$ has height $1$ in
$A / \mfrak{p}_{i}$.

\begin{lem} \label{lem10.1}
If $\xi= (\mfrak{p}, \ldots )$ is a saturated chain then
$A_{\xi}$ is a Beilinson completion algebra.
\end{lem}

\begin{proof}
By \cite{Ye1} Corollary 3.3.5,
$A_{\xi}$ is a complete semi-local noetherian ring with
Jacobson radical $\mfrak{p}_{\xi}$, and
$A_{\xi} / \mfrak{p}_{\xi} = K_{\xi}$,
where $K := A_{\mfrak{p}} / \mfrak{p}_{\mfrak{p}}$.
Choose a coefficient field
$\sigma : K \ar \widehat{A}_{\mfrak{p}} = A_{(\mfrak{p})}$.
By \cite{Ye1} Proposition 3.3.6, $K_{\xi}$ is a finite
product of TLFs, and by the proof of \cite{Ye1} Theorem 3.3.8,
$\sigma$ extends to a lifting
$\sigma_{\xi} : K_{\xi} \ar A_{\xi}$. Sending
$t_{1}, \ldots, t_{m}$ to generators of the ideal
$\mfrak{p}$, we get a topologically strict surjection
$K_{\xi} [[ t_{1}, \ldots, t_{m} ]] \surj A_{\xi}$.
\end{proof}

\section{Construction of the Residue Complex $\mcal{K}^{\bdot}_{X}$}

Let $k$ be a perfect field, and let $X$ be a scheme of finite type
over $k$. By a {\em chain of points} in $X$ we mean a sequence
$\xi = (x_{0}, \ldots, x_{n})$ of points with
$x_{i + 1} \in \overline{\{x_{i}\}}$.

\begin{dfn}
For any quasi-coherent $\mcal{O}_{X}$-module $\mcal{M}$,
define the Beilinson completion $\mcal{M}_{\xi}$
by taking an affine open neighborhood $U = \opn{Spec} A \subset X$
of $x_{n}$, and setting
$\mcal{M}_{\xi} := \Gamma(U, \mcal{M})_{\xi}$
as in Definition \ref{dfn10.1}.
\end{dfn}

These completions first appeared as the local factors of Beilinson's
adeles in \cite{Be}, and were studied in detail in \cite{Ye1}.

According to Lemma \ref{lem10.1}, if
$\xi = (x_{0}, \ldots, x_{n})$ is saturated, i.e.\
$\overline{\{x_{i + 1}\}} \subset \overline{\{x_{i}\}}$
has codimension $1$, then
$\mcal{O}_{X, \xi}$ is a BCA.
We shall be interested in the covertex maps
\[ \begin{array}{rcl}
\partial^{-} :  &  \mcal{O}_{X,(x_{0})} & \ar \mcal{O}_{X,\xi}
\\
\partial^{+} :  &  \mcal{O}_{X,(x_{n})} & \ar \mcal{O}_{X,\xi}
\end{array} \]
which arise naturally from the completion process
(cf.\ \cite{Ye1} \S 3.1).

\begin{lem} \label{lem1.1}
$\partial^{+}$ is flat, topologically \'{e}tale relative to $k$, and a
morphism in $\msf{BCA}(k)$.
$\partial^{-}$ is an intensification homomorphism.
\end{lem}

\begin{proof}
By definition
$\partial^{-} = \partial^{1} \circ \cdots \circ \partial^{n}$ and
$\partial^{+} = \partial^{0} \circ \cdots \circ \partial^{0}$, where
$\partial^{i} : \mcal{O}_{X,\partial_{i} \xi} \ar \mcal{O}_{X,\xi}$
is the $i$-th coface operator.
First let us prove that
$\partial^{0} : \mcal{O}_{X,\partial_{0} \xi} \ar \mcal{O}_{X,\xi}$
is a morphism of BCAs. This follows from \cite{Ye1} Theorem 3.3.2 (d),
since we may assume that $X$ is integral with generic point $x_{0}$.
By part (b) of the same theorem,
$\partial^{n} : \mcal{O}_{X,\partial_{n} \xi} \ar
\mcal{O}_{X, \xi}$
is finitely ramified and radically unramified (in the sense of \cite{Ye2}
Definition 3.1).

Now according to \cite{Ye1} Corollary 3.2.8,
$\partial^{i} :  \mcal{O}_{X,\partial_{i} \xi} \ar \mcal{O}_{X,\xi}$
is topologically \'{e}tale relative to $k$, for any $i$. We claim
it is also flat. For $i=0$,
$\mcal{O}_{X,\partial_{0} \xi} \ar
(\mcal{O}_{X,\partial_{0} \xi})_{x_{0}}
=  (\mcal{O}_{X, x_{0}})_{\partial_{0} \xi}$
is a localization, so it's flat. The map from
$(\mcal{O}_{X, x_{0}})_{\partial_{0} \xi}$ to its
$\mfrak{m}_{x_{0}}$-adic
completion $\mcal{O}_{X, \xi}$ is also flat (these rings are
noetherian). For $i > 0$, by induction on the length of chains,
$\mcal{O}_{X,\partial_{0} \partial_{i} \xi} \ar
\mcal{O}_{X,\partial_{0} \xi}$
is flat, and hence so is
$(\mcal{O}_{X, x_{0}})_{\partial_{0} \partial_{i} \xi} \ar
(\mcal{O}_{X, x_{0}})_{\partial_{0} \xi}$. Now use
\cite{CA} Chapter III \S 5.4 Proposition 4 to conclude that
\[ \mcal{O}_{X,\partial_{i} \xi} =
\lim_{\leftarrow j} (\mcal{O}_{X, x_{0}} / \mfrak{m}_{x_{0}}^{j})_{
\partial_{0} \partial_{i} \xi} \ar
\lim_{\leftarrow j} (\mcal{O}_{X, x_{0}} / \mfrak{m}_{x_{0}}^{j})_{
\partial_{0} \xi} = \mcal{O}_{X, \xi} \]
 is flat.
\end{proof}

\begin{exa} \label{1.1}
Take $X = \mbf{A}^{2} := \opn{Spec} k[s,t]$,
$x := (0)$,
$y := (t)$ and $z := (s,t)$. Then with
$L := k(x) = \mcal{O}_{X, (x)}$,
$\widehat{L} := k(x)_{(y)} = \mcal{O}_{X, (x, y)}$,
$A := \mcal{O}_{X, (y)}$,
$\widehat{A} := \mcal{O}_{X, (y, z)}$
$K := k(y)$ and
$\widehat{K} := k(y)_{(z)}$
we are exactly in the situation of Example \ref{exa10.1}.
\end{exa}

\begin{dfn}
Given  a point $x$ in $X$, let
$\mcal{K}_{X}(x)$ be the skyscraper sheaf with support
$\overline{\{ x \}}$ and group of sections
$\mcal{K}(\mcal{O}_{X,(x)})$.
\end{dfn}

The sheaf $\mcal{K}_{X}(x)$ is a quasi-coherent
$\mcal{O}_{X}$-module,
and is an injective hull of $k(x)$ in the category $\msf{Mod}(X)$ of
$\mcal{O}_{X}$-modules.

\begin{dfn} \label{dfn1.1}
Given a saturated chain $\xi = (x, \ldots, y)$ in $X$, define an
$\mcal{O}_{X}$-linear homomorphism
$\delta_{\xi} : \mcal{K}_{X}(x) \ar \mcal{K}_{X}(y)$, called the
coboundary map along $\xi$, by
\[ \delta_{\xi} :  \mcal{K}(\mcal{O}_{X,(x)})
\xrightarrow{ \mrm{q}_{\partial^{-}} }
\mcal{K}(\mcal{O}_{X,\xi}) \xrightarrow{ \opn{Tr}_{\partial^{+}} }
\mcal{K}(\mcal{O}_{X,(y)}) . \]
\end{dfn}

Throughout sections 2 and 3 the following convention shall be used.
Let $f: X \ar Y$ be a morphism of schemes, and let $x \in X$,
$y \in Y$ be points. We will
write $x \mid y$ to indicate that $x$ is a closed point in the fiber
$X_{y} := X \times_{Y} \opn{Spec} k(y) \cong f^{-1}(y)$. Similarly for
chains: we write
$(x_{0}, \ldots, x_{n}) \mid (y_{0}, \ldots, y_{n})$ if
$x_{i} \mid y_{i}$ for all $i$.

\begin{lem}
Suppose $x \mid y$. Then
$f^{*} : \mcal{O}_{Y,(y)} \ar \mcal{O}_{X,(x)}$ is a
morphism of BCAs. If $f$ is quasi-finite then $f^{*}$ is finite, and
if $f$ is \'{e}tale then $f^{*}$ is an intensification.
\end{lem}

\begin{proof}
Immediate from the definitions.
\end{proof}

\begin{lem} \label{lem1.2}
Suppose $f: X \ar Y$ is a quasi-finite morphism. Let
$\eta = (y_{0}, \ldots, y_{n})$ be a saturated chain in $Y$ and
let $x \in X$, $x \mid y_{n}$. Consider the \textup{(}finite\textup{)}
set of chains in $X$:
\[ \Xi := \{ \xi = (x_{0}, \ldots, x_{n})\  \mid  \
\xi \mid \eta \textup{ and } x_{n} = x \}. \]
Then there is a canonical isomorphism  of BCAs
\[ \prod_{\xi \in \Xi} \mcal{O}_{X,\xi} \cong
\mcal{O}_{Y,\eta} \otimes_{\mcal{O}_{Y,(y_{n})}}
\mcal{O}_{X,(x)} . \]
\end{lem}

\begin{proof}
Set
$\widehat{X} := \opn{Spec} \mcal{O}_{X,(x)}$ and
$\widehat{Y} := \opn{Spec} \mcal{O}_{Y,(y_{n})}$, so the induced
morphism
$\widehat{f} : \widehat{X} \ar \widehat{Y}$ is finite.
Let us denote by $\xi, \widehat{\xi}, \widehat{\eta}$ variable saturated
chains in $X, \widehat{X}, \widehat{Y}$ respectively. For any
$\widehat{\eta} \mid \eta$ one has
\begin{equation} \label{eqn1.1}
\prod_{\widehat{\xi} \mid \widehat{\eta}} \mcal{O}_{\widehat{X},
\widehat{\xi}} \cong
\mcal{O}_{\widehat{Y}, \widehat{\eta}}
\otimes_{\mcal{O}_{Y, (y_{n})}} \mcal{O}_{X, (x)},
\end{equation}
by \cite{Ye1} Proposition 3.2.3; note that the completion is defined
on any noetherian scheme. Now from ibid.\ Corollary 3.3.13 one has
$\mcal{O}_{X, \xi} \cong \prod_{\widehat{\xi} \mid \xi}
\mcal{O}_{\widehat{X}, \widehat{\xi}}$, so taking the product over all
$\xi \in \Xi$ and $\widehat{\eta} \mid \eta$ the lemma is proved.
\end{proof}

\begin{dfn} \label{dfn1.2}
Given an \'{e}tale morphism $g : X \ar Y$ and a point $x \in X$, let
$y := g(x)$, so
$g^{*} : \mcal{O}_{Y, (y)} \ar \mcal{O}_{X, (x)}$
is an intensification. Define
\[ \mrm{q}_{g} : \mcal{K}_{Y}(y) \ar g_{*} \mcal{K}_{X}(x) \]
to be the $\mcal{O}_{Y}$-linear homomorphism corresponding to
$\mrm{q}_{g^{*}} : \mcal{K}(\mcal{O}_{Y, (y)}) \ar
\mcal{K}(\mcal{O}_{X, (x)})$
of \cite{Ye2} Proposition 7.2.
\end{dfn}

\begin{prop} \label{prop1.1}
Let $g : X \ar Y$ be an \'{e}tale morphism.
\begin{enumerate}
\rmitem{a} Given a point $y \in Y$, the homomorphism
$1 \otimes \mrm{q}_{g} : g^{*} \mcal{K}_{Y}(y) \ar$ \blnk{3mm} \linebreak
$\bigoplus_{x \mid y}
\mcal{K}_{X}(x)$ is an isomorphism of $\mcal{O}_{X}$-modules.

\rmitem{b} Let
$\eta = (y_{0}, \ldots, y_{n})$ be a saturated chain in $Y$. Then
\[ (1 \otimes \mrm{q}_{g}) \circ g^{*}(\delta_{\eta}) =
(\sum_{\xi \mid \eta} \delta_{\xi}) \circ (1 \otimes \mrm{q}_{g}) \]
as homomorphisms
$g^{*} \mcal{K}_{Y}(y_{0}) \ar \bigoplus_{x_{n} \mid y_{n}}
\mcal{K}_{X}(x_{n})$.
\end{enumerate}
\end{prop}

\begin{proof}
(a)\ Because $\mcal{K}_{Y}(y)$ is an artinian
$\mcal{O}_{Y,y}$-module,
and $g$ is quasi-finite, we find that
\[ g^{*} \mcal{K}_{Y}(y) = \bigoplus_{x \mid y} \mcal{O}_{X, (x)}
\otimes_{\mcal{O}_{Y, (y)}} \mcal{K}_{Y}(y)\ . \]
Now use \cite{Ye2} Proposition 7.2 (i).

\medskip \noindent (b)\
From Lemma \ref{lem1.2} and from \cite{Ye2} Theorem 3.8 we see
that for every $\xi \mid \eta$,
$f^{*} : \mcal{O}_{Y,\eta} \ar \mcal{O}_{X, \xi}$ is both a finite
morphism and an intensification.
By the definition of the coboundary maps, it suffices to verify that the
diagram

\bigskip \noindent
\[ \begin{CD}
\mcal{K}(\mcal{O}_{Y,(y_{0})})  @>{\mrm{q}}>>
\mcal{K}(\mcal{O}_{Y,\eta}) @>{\opn{Tr}}>>
\mcal{K}(\mcal{O}_{Y,(y_{n})}) \\
@V{\mrm{q}}VV @V{\mrm{q}}VV @V{\mrm{q}}VV \\
\bigoplus_{x_{0} \mid y_{0}} \mcal{K}(\mcal{O}_{X,(x_{0})})
@>{\mrm{q}}>> \bigoplus_{\xi \mid \eta} \mcal{K}(\mcal{O}_{X,\xi})
@>{\opn{Tr}}>>
\bigoplus_{x_{n} \mid y_{n}} \mcal{K}(\mcal{O}_{X,(x_{n})})
\end{CD} \]

\medskip \noindent
commutes. The left square commutes by \cite{Ye2} Proposition 7.2 (iv),
whereas the right square commutes by Lemma \ref{lem1.2} and
\cite{Ye2} Theorem 7.4 (ii).
\end{proof}

\begin{dfn} \label{dfn1.3}
Let $f: X \ar Y$ be a morphism between finite type $k$-schemes,
let $x \in X$ be a point, and let $y := f(x)$. Define an
$\mcal{O}_{Y}$-linear homomorphism
\[ \opn{Tr}_{f} : f_{*} \mcal{K}_{X}(x) \ar \mcal{K}_{Y}(y) \]
as follows:
\begin{enumerate}
\rmitem{i} If $x$ is closed in its fiber $X_{y}$, then
$f^{*} : \mcal{O}_{Y,(y)} \ar \mcal{O}_{X,(x)}$ is a morphism in
$\msf{BCA}(k)$. Let $\opn{Tr}_{f}$
be the homomorphism corresponding to
$\opn{Tr}_{f^{*}} : \mcal{K}(\mcal{O}_{X,(x)}) \ar
\mcal{K}(\mcal{O}_{Y,(y)})$
of \cite{Ye2} Theorem 7.4.

\rmitem{ii} If $x$ is not closed in its fiber, set $\opn{Tr}_{f} := 0$.
\end{enumerate}
\end{dfn}

\begin{prop} \label{prop1.2}
Let $f : X \ar Y$ be a finite morphism.
\begin{enumerate}
\rmitem{a} For any $y \in Y$ the homomorphism of
$\mcal{O}_{Y}$-modules
\[ \bigoplus_{x \mid y} f_{*} \mcal{K}_{X}(x) \ar
\mcal{H}om_{\mcal{O}_{Y}}(f_{*} \mcal{O}_{X},
\mcal{K}_{Y}(y)) \]
induced by $\opn{Tr}_{f}$ is an isomorphism.

\rmitem{b} Let $\eta = (y_{0}, \ldots, y_{n})$ be a saturated
chain in $Y$. Then
\[ \delta_{\eta} \opn{Tr}_{f} =
\opn{Tr}_{f} \sum_{\xi \mid \eta} f_{*}(\delta_{\xi}) :
\bigoplus_{x_{0} \mid y_{0}} f_{*} \mcal{K}_{X}(x_{0}) \ar
\mcal{K}_{Y}(y_{n}) . \]
The sums are over saturated chains $\xi = (x_{0}, \ldots, x_{n})$
in $X$.
\end{enumerate}
\end{prop}

\begin{proof}
(a)\
By \cite{Ye1} Proposition 3.2.3 we have
$\prod_{x \mid y} \mcal{O}_{X, (x)} \cong
(f_{*} \mcal{O}_{X})_{(y)}$.
Now use \cite{Ye2} Theorem 7.4 (iv).

\medskip \noindent (b)\
Use the same diagram which appears in the proof of Proposition
\ref{prop1.1},
only reverse the vertical arrows and label them $\opn{Tr}_{f}$. Then the
commutativity follows from \cite{Ye2} Theorem 7.4  (i), (ii).
\end{proof}

In \cite{Ye1} \S 4.3 the notion of a system of residue data on a reduced
scheme was introduced.

\begin{prop} \label{prop1.3}
Suppose $X$ is a reduced scheme. Then
$(\{ \mcal{K}_{X}(x) \}, \{ \delta_{\xi} \},$ \linebreak
$\{ \Psi_{\sigma}^{-1} \})$,
where $x$ runs over the points of $X$, $\xi$ runs over the
saturated chains
in $X$, and $\sigma : k(x) \ar \mcal{O}_{X,(x)}$ runs over all
possible coefficient fields, is a system of residue data on $X$.
\end{prop}

\begin{proof}
We must check condition (\dag) of \cite{Ye1} Definition 4.3.10.
So let $\xi = (x, \ldots, y)$ be a saturated chain, and let
$\sigma : k(x) \ar \mcal{O}_{X,(x)}$ and
$\tau : k(y) \ar \mcal{O}_{X,(y)}$ be compatible coefficient fields.
Denote also by $\tau$ the composed morphism
$\partial^{+} \tau : k(y) \ar \mcal{O}_{X,\xi}$.
Then we get a coefficient field
$\sigma_{\xi} : k(\xi) = k(x)_{\xi} \ar \mcal{O}_{X,\xi}$
extending $\sigma$, which is a $k(y)$-algebra map via $\tau$.
Consider the diagram:

\bigskip \noindent
\[ \begin{CD}
\mcal{K}(\mcal{O}_{X,(x)}) @>\mrm{q}>>
\mcal{K}(\mcal{O}_{X,\xi}) @>=>>
\mcal{K}(\mcal{O}_{X,\xi}) @>\opn{Tr}>>
\mcal{K}(\mcal{O}_{X,(y)}) \\
@V{\Psi_{\sigma}}VV  @V{\Psi_{\sigma_{\xi}}}VV
@V{\Psi_{\tau}}VV @V{\Psi_{\tau}}VV \\
\opn{Dual}_{\sigma} \mcal{O}_{X,(x)}    @>{\mrm{q}_{\sigma}}>>
\opn{Dual}_{\sigma_{\xi}} \mcal{O}_{X,\xi} @>{h}>>
\opn{Dual}_{\tau} \mcal{O}_{X,\xi}
@>{\opn{Tr}_{\tau}}>>
\opn{Dual}_{\tau} \mcal{O}_{X, (y)}
\end{CD} \]

\medskip \noindent
where for a $k(\xi)$-linear homomorphism
$\phi : \mcal{O}_{X,\xi} \ar \omega(k(\xi))$,
$h(\phi) = $ \blnk{3mm} \linebreak
$\opn{Res}_{k(\xi) / k(y)} \phi$
(cf.\ \cite{Ye2} Theorem 6.14).
The left square commutes by \cite{Ye2} Proposition 7.2 (iii); the
middle square
commutes by ibid.\ Theorem 6.14 (i); and the right square commutes
by ibid.\ Theorem 7.4  (i), (iii). But going along the bottom
of the diagram we get
$\opn{Tr}_{\tau} h \mrm{q}_{\sigma} = \delta_{\xi, \sigma / \tau}$,
as defined in \cite{Ye1} Lemma 4.3.3.
\end{proof}

\begin{lem} \label{lem1.5}
Let $\xi = (x, \ldots, y)$ and $\eta = (y, \ldots, z)$  be saturated
chains in the scheme $X$, and let
$\xi \vee \partial_{0} \eta := (x, \ldots, y, \ldots, z)$ be the
concatenated chain. Then
there is a canonical isomorphism of BCAs
\[ \mcal{O}_{X,\xi \vee \partial_{0} \eta} \cong \mcal{O}_{X,\xi}
\otimes^{(\wedge)}_{\mcal{O}_{X,(y)}} \mcal{O}_{X,\eta} \]
\textup{(}intensification base change\textup{)}.
\end{lem}

\begin{proof}
Choose a coefficient field $\sigma: k(y) \ar \mcal{O}_{X,(y)}$. This
induces a coefficient ring
$\sigma_{\eta}: k(\eta) \ar \mcal{O}_{X,\eta}$, and using \cite{Ye2}
Theorem 3.8 and \cite{Ye1} Theorem 4.1.12 one gets
\[ \mcal{O}_{X,\xi \vee \partial_{0} \eta} \cong
\mcal{O}_{X,\xi} \otimes^{(\wedge)}_{k(y)} k(\eta) \cong
\mcal{O}_{X,\xi} \otimes^{(\wedge)}_{\mcal{O}_{X,(y)}}
\mcal{O}_{X,\eta} . \]
\end{proof}

\begin{lem} \label{lem1.3}
\begin{enumerate}
\item Let $\xi = (x, \ldots, y)$ and
$\eta = (y, \ldots, z)$ be saturated chains and let
$\xi \vee \partial_{0} \eta = (x, \ldots, y, \ldots, z)$. Then
$\delta_{\eta} \delta_{\xi} = \delta_{\xi \vee \partial_{0} \eta}$.

\item Given a point $x \in X$ and an element
$\alpha \in \mcal{K}_{X}(x)$,
for all but finitely
many saturated chains $\xi = (x, \ldots)$ in $X$ one has
$\delta_{\xi}(\alpha) = 0$.

\item \textup{(}Residue Theorem\textup{)} Let $x,z \in X$ be points
s.t.\ $z \in \overline{\{ x \}}$ and
$\opn{codim}(\overline{\{ z \}},$ \linebreak
$\overline{\{ x \}}) = 2$.
Then
$\sum_{y}\  \delta_{(x,y,z)} = 0$.
\end{enumerate}
\end{lem}

\begin{proof}
Using Lemma \ref{lem1.5} we see that part
1 is a consequence of the base change property of traces,  cf.\
\cite{Ye2} Theorem 7.4 (ii). Assertions
2 and 3 are local, by Proposition \ref{prop1.1}, so we may assume
there is a closed immersion
$f : X \ar \mbf{A}_{k}^{n}$. By Proposition
\ref{prop1.2}, we can replace $X$ with $\mbf{A}_{k}^{n}$, and
so assume
that $X$ is reduced. Now according to Proposition  \ref{prop1.3}
and \cite{Ye1} Lemma 4.3.19, both 2 and 3 hold.
\end{proof}

\begin{dfn}
For any integer $q$ define a quasi-coherent sheaf
\[  \mcal{K}_{X}^{q} :=
\bigoplus_{\opn{dim} \overline{\{x\}} = -q} \mcal{K}_{X}(x) . \]
By Lemma \ref{lem1.3} there an $\mcal{O}_{X}$-linear homomorphism
\[ \delta := (-1)^{q+1} \sum_{(x,y)} \delta_{(x,y)} :
\mcal{K}_{X}^{q} \ar \mcal{K}_{X}^{q+1} , \]
satisfying  $\delta^{2} = 0$. The complex
$(\mcal{K}_{X}^{\bdot}, \delta)$ is called the {\em Grothendieck
residue complex} of $X$.
\end{dfn}

In Corollary \ref{cor2.2} we will prove that $\mcal{K}_{X}^{\bdot}$
is canonically isomorphic (in the derived category $\msf{D}(X)$)
to $\pi^{!} k$, where $\pi : X \ar \opn{Spec} k$ is the structural
morphism.

\begin{rem}
A heuristic for the negative grading of $\mcal{K}_{X}^{\bdot}$
and the sign $(-1)^{q+1}$ is that the residue complex is the
``$k$-linear dual'' of a hypothetical ``complex of localizations''
$\cdots \prod \mcal{O}_{X,y} \ar \prod \mcal{O}_{X,x} \ar \cdots$.
Actually, there is a naturally defined complex which is built up from
all localizations and completions: the Beilinson adeles
$\ul{\mbb{A}}^{\bdot}_{\mrm{red}}(\mcal{O}_{X})$
(cf.\ \cite{Be} and \cite{HY1}).
$\ul{\mbb{A}}^{\bdot}_{\mrm{red}}(\mcal{O}_{X})$
is a DGA, and $\mcal{K}_{X}^{\bdot}$ is naturally a right DG-module
over it. See \cite{Ye4}, and cf.\ also Remark \ref{rem5.3}.
\end{rem}

\begin{dfn} \label{eqn1.4}
\begin{enumerate}
\item Let $f:X \ar Y$ be a morphism of schemes. Define a homomorphism of
graded $\mcal{O}_{Y}$-modules
$\opn{Tr}_{f} : f_{*} \mcal{K}_{X}^{\bdot} \ar \mcal{K}_{Y}^{\bdot}$
by summing the local trace maps of Definition \ref{dfn1.3}.
\item Let $g : U \ar X$ be an \'{e}tale morphism. Define
$\mrm{q}_{g} : \mcal{K}_{X}^{\bdot} \ar g_{*} \mcal{K}_{U}^{\bdot}$
by summing all local homomorphisms $\mrm{q}_{g}$
of Definition \ref{dfn1.2}.
\end{enumerate}
\end{dfn}

\begin{thm} \label{thm1.1}
Let $X$ be a $k$-scheme of finite type.
\begin{enumerate}
\rmitem{a} $(\mcal{K}_{X}^{\bdot}, \delta)$ is a residual complex
on $X$ \textup{(}cf.\ \cite{RD} Chapter \textup{VI \S 1)}.

\rmitem{b} If $g : U \ar X$ is an \'{e}tale morphism, then
$1 \otimes \mrm{q}_{g} : g^{*} \mcal{K}_{X}^{\bdot} \ar
\mcal{K}_{U}^{\bdot}$
is an isomorphism of complexes.

\rmitem{c} If $f : X \ar Y$ is a finite morphism, then
$\opn{Tr}_{f} : f_{*} \mcal{K}_{X}^{\bdot} \ar \mcal{K}_{Y}^{\bdot}$
is a homomorphism of complexes, and the induced map
\[ f_{*} \mcal{K}_{X}^{\bdot} \ar
\mcal{H}om_{\mcal{O}_{Y}}(f_{*} \mcal{O}_{X},
\mcal{K}_{Y}^{\bdot}) \]
is an isomorphism of complexes.

\rmitem{d} If $X$ is reduced, then $(\mcal{K}_{X}^{\bdot}, \delta)$
is canonically isomorphic to the complex constructed in \cite{Ye1} \S
\textup{4.3}.
In particular, if $X$ is smooth irreducible of dimension $n$, there
is a quasi-isomorphism
$\mrm{C}_{X} : \Omega^{n}_{X/k}[n] \ar \mcal{K}^{\bdot}_{X}$.
\end{enumerate}
\end{thm}

\begin{proof}
Parts (b), (c), (d) are immediate consequences of Propositions
\ref{prop1.1}, \ref{prop1.2} and \ref{prop1.3} here, and \cite{Ye1}
Theorem 4.5.2. (Note that the sign of $\delta$ in \cite{Ye1} is
different.) As for part (a),
clearly $\mcal{K}_{X}^{\bdot}$ is a direct sum of injective hulls of
all the residue fields in $X$, with multiplicities $1$. It remains to
prove that $\mcal{K}_{X}^{\bdot}$ has coherent cohomology.
Since this is a local
question, we can assume using part (b) that $X$ is a closed subscheme of
$\mbf{A}_{k}^{n}$. According to parts (c) and (d) of this theorem and
\cite{Ye1} Corollary 4.5.6, $\mcal{K}_{X}^{\bdot}$ has
coherent cohomology.
\end{proof}

From part (b) of the theorem we get:

\begin{cor}
The presheaf
$U \mapsto \Gamma(U, \mcal{K}_{U}^{\bdot})$
is a sheaf on $X_{\mrm{et}}$, the small \'{e}tale
site over $X$.
\end{cor}

\begin{dfn} \label{dfn1.6}
For an $\mcal{O}_{X}$-module $\mcal{M}$
define dual complex
\[ \opn{Dual}_{X} \mcal{M} :=
\mcal{H}om^{\bdot}_{\mcal{O}_{X}}(\mcal{M}, \mcal{K}_{X}^{\bdot}) . \]
\end{dfn}

Observe that since $\mcal{K}_{X}^{\bdot}$ is complex of injectives the
derived functor
$\opn{Dual}_{X} : \msf{D}(X)^{\circ} \ar \msf{D}(X)$
is defined. Moreover, since $\mcal{K}_{X}^{\bdot}$ is dualizing, the
adjunction morphism
$1 \ar \opn{Dual}_{X} \opn{Dual}_{X}$
is an isomorphism on
$\msf{D}_{\mrm{c}}^{\mrm{b}}(X)$.
We shall sometimes write $\opn{Dual} \mcal{M}$ instead of
$\opn{Dual}_{X} \mcal{M}$.

\section{Duality for Proper Morphisms}

In this section we prove that if $f : X \ar Y$ is a proper morphism of
$k$-schemes, then the trace map $\opn{Tr}_{f}$
of Definition \ref{dfn1.3} is a homomorphism of complexes, and it induces a
duality in the derived categories.

\begin{prop} \label{prop2.1}
Let $f: X \ar Y$ be a proper morphism between finite type $k$-schemes, and
let $\eta = (y_{0}, \ldots, y_{n})$ be  a saturated chain in $Y$. Then there
exists a canonical isomorphism of BCAs
\[ \prod_{\xi \mid \eta} \mcal{O}_{X,\xi} \cong
\prod_{x_{0} \mid y_{0}} \mcal{O}_{X,(x_{0})} \otimes^{(\wedge)}_{
\mcal{O}_{Y,(y_{0})}} \mcal{O}_{Y,\eta}\ , \]
where $\xi = (x_{0}, \ldots, x_{n})$ denotes a variable chain in $X$
lying over $\eta$.
\end{prop}

\begin{proof}
The proof is by induction on $n$. For $n=0$ this is trivial. Assume $n=1$.
Let
$Z := \overline{\{ x_{0} \}}_{\mrm{red}}$,
so $\mcal{O}_{Z, x_{1}}$ is a $1$-dimensional local ring inside
$k(Z) = k(x_{0})$. Considering the integral closure of
$\mcal{O}_{Z, x_{1}}$ we see that
$k((x_{0}, x_{1})) = k(x_{0})_{(x_{1})} =
k(x_{0}) \otimes \mcal{O}_{Z, (x_{1})}$
is the product of the completions of $k(x_{0})$ at all discrete valuations
centered on $x_{1} \in Z$ (cf.\ \cite{Ye1} Theorem 3.3.2).
So by the valuative criterion for properness we get
\begin{equation} \label{eqn2.3}
\prod_{(x_{0}, x_{1}) \mid (y_{0}, y_{1})} k(x_{0})_{(x_{1})} \cong
\prod_{x_{0} \mid y_{0}} k(x_{0}) \otimes_{k(y_{0})} k(y_{0})_{(y_{1})}.
\end{equation}
For $i \geq 1$ the  morphism of BCAs
\[ \prod_{x_{0} \mid y_{0}} (\mcal{O}_{X, x_{0}} /
\mfrak{m}_{x_{0}}^{i})
\otimes_{\mcal{O}_{Y, (y_{0})}} \mcal{O}_{Y, (y_{0}, y_{1})} \ar
\prod_{(x_{0}, x_{1}) \mid (y_{0}, y_{1})}
\mcal{O}_{X, (x_{0}, x_{1})} / \mfrak{m}_{(x_{0}, x_{1})}^{i} \]
is bijective, since both sides are flat over
$\mcal{O}_{X, x_{0}} / \mfrak{m}_{x_{0}}^{i}$, and by equation
(\ref{eqn2.3}) (cf.\ \cite{Ye2} Proposition 3.5).
Passing to the inverse limit in $i$ we get an isomorphism of BCAs
\begin{equation} \label{eqn2.4}
\prod_{x_{0} \mid y_{0}} \mcal{O}_{X, (x_{0})} \otimes^{(\wedge)}_{
\mcal{O}_{Y,(y_{0})}} \mcal{O}_{Y,(y_{0}, y_{1})} \cong
\prod_{(x_{0}, x_{1}) \mid (y_{0}, y_{1})} \mcal{O}_{X,(x_{0}, x_{1})}\ .
\end{equation}

Now suppose $n \geq 2$. Then we get
\[ \begin{array}{l}
\prod_{x_{0} \mid y_{0}} \mcal{O}_{X,(x_{0})} \otimes^{(\wedge)}_{
\mcal{O}_{Y,(y_{0})}} \mcal{O}_{Y, \eta} \\[2mm]
\begin{array}{rclc}
\blnk{25mm} & \cong &
\prod_{(x_{0}, x_{1}) \mid (y_{0}, y_{1})} \mcal{O}_{X,(x_{0}, x_{1})}
\otimes^{(\wedge)}_{\mcal{O}_{Y,(y_{0}, y_{1})}} \mcal{O}_{Y, \eta}
& \text{(i)} \\[2mm]
& \cong &
\prod_{(x_{0}, x_{1}) \mid (y_{0}, y_{1})} \mcal{O}_{X,(x_{0}, x_{1})}
\otimes^{(\wedge)}_{\mcal{O}_{Y,(y_{1})}}
\mcal{O}_{Y, \partial_{0} \eta}
& \text{(ii)} \\[2mm]
& \cong &
\prod_{\xi \mid \eta} \mcal{O}_{X,(x_{0}, x_{1})}
\otimes^{(\wedge)}_{\mcal{O}_{X,(x_{1})}}
\mcal{O}_{X, \partial_{0} \xi}
& \text{(iii)} \\[2mm]
& \cong &
\prod_{\xi \mid \eta} \mcal{O}_{X, \xi}
& \text{(iv)}
\end{array}
\end{array} \]
where associativity of intensification base change (\cite{Ye2}
Proposition 3.10)
is used repeatedly; in (i) we use formula (\ref{eqn2.4}); in (ii) we use
Lemma \ref{lem1.5} applied to $\mcal{O}_{Y, \eta}$; in (iii) we use
the induction hypothesis; and (iv) is another application of
Lemma \ref{lem1.5}.
\end{proof}

The next theorem is our version of \cite{RD} Ch.\ VII Theorem 2.1:

\begin{thm} \label{thm2.1}
\textup{(Global Residue Theorem)}\
Let $f : X \ar Y$ be a proper morphism between $k$-schemes of finite type.
Then
$\opn{Tr}_{f}: f_{*} \mcal{K}_{X}^{\bdot} \ar \mcal{K}_{Y}^{\bdot}$
is a homomorphism of complexes.
\end{thm}

\begin{proof}
Fix a point $x_{0} \in X$, and let $y_{0} := f(x_{0})$. First assume that
$x_{0}$ is closed in its fiber $X_{y_{0}} = f^{-1}(y_{0})$.
Let $y_{1}$ be an immediate specialization of $y_{0}$.
By Proposition \ref{prop2.1} we have
\[ \prod_{x_{1} \mid y_{1}} \mcal{O}_{X, (x_{0}, x_{1})} \cong
\mcal{O}_{X, (x_{0})} \otimes^{(\wedge)}_{\mcal{O}_{Y, (y_{0})}}
\mcal{O}_{Y, (y_{0}, y_{1})}, \]
so just as in Proposition \ref{prop1.2} (b), we get
\[ \delta_{(y_{0}, y_{1})}  \opn{Tr}_{f} =
\sum_{x_{1} \mid y_{1}} \opn{Tr}_{f} \delta_{(x_{0}, x_{1})} :
f_{*} \mcal{K}_{X}(x_{0}) \ar \mcal{K}_{Y}(y_{1})  . \]

Next assume  $x_{0}, y_{0}$ are as above, but $x_{0}$ is not closed in
the fiber $X_{y_{0}}$. The only possibility to have an immediate
specialization $x_{1}$ of $x_{0}$ which is closed in its fiber, is if
$x_{1} \in X_{y_{0}}$ and
$Z := \overline{\{ x_{0} \}}_{\mrm{red}} \subset X_{y_{0}}$
is a curve. We have to show that
\begin{equation} \label{eqn2.5}
\sum_{x_{1} \mid y_{0}} \opn{Tr}_{f} \delta_{(x_{0}, x_{1})} = 0
: f_{*} \mcal{K}_{X}(x_{0}) \ar \mcal{K}_{Y}(y_{0})  .
\end{equation}

Since
$\mcal{K}_{Z}(x_{0}) \subset \mcal{K}_{X}(x_{0})$
is an essential submodule over $\mcal{O}_{Y, y_{0}}$ it suffices to
check (\ref{eqn2.5}) on $\mcal{K}_{Z}(x_{0})$. Thus we may assume
$X = \overline{\{ x_{0} \}}_{\mrm{red}}$
and
$Y = \overline{\{ y_{0} \}}_{\mrm{red}}$.
After factoring $X \ar Y$ through a suitable finite radiciel morphism
$X \ar \tilde{X}$, and using Proposition \ref{prop1.2}, we may further
assume that $K = k(Y) \ar k(X)$ is separable. Now
\[  \opn{Tr}_{f} \delta_{(x_{0}, x_{1})} =
\opn{Res}_{k((x_{0}, x_{1})) / K} : \Omega^{n+1}_{k(X) / k} \ar
\Omega^{n}_{K / k} . \]
Since
$\Omega^{n+1}_{k(X) / k} = \Omega^{1}_{k(X) / k} \wedge
\Omega^{n}_{K / k}$,
it suffices to check that
\[ \sum_{x_{1} \in X} \opn{Res}_{k((x_{0}, x_{1})) / K} = 0 :
\Omega^{1}_{k(X) / k} \ar K . \]

Let $K'$ be the maximal purely inseparable extension of $K$ in
an algebraic closure, and let
$X' := X \times_{K} K'$. So
\[ k((x_{0}, x_{1})) \otimes_{K} K' \cong
\prod_{(x'_{0}, x'_{1})  \mid (x_{0}, x_{1})} k((x'_{0}, x'_{1})) \]
where $(x'_{0}, x'_{1})$ are chains in $X'$.
According to \cite{Ye1} Lemma 2.4.14 we may assume $k = K = K'$.
Since now $K$ is perfect, we are in the position to use the
well known Residue Theorem for curves (cf.\ \cite{Ye1} Theorem 4.2.15).
\end{proof}

\begin{cor} \label{cor2.1}
Let $f : X \ar Y$ be a morphism between $k$-schemes of finite type, and
let $Z \subset X$ be a closed subscheme which is proper over $Y$. Then
$\opn{Tr}_{f}: f_{*} \Gamma_{Z} \mcal{K}_{X}^{\bdot} \ar
\mcal{K}_{Y}^{\bdot}$
is a homomorphism of complexes.
\end{cor}

\begin{proof}
Let $\mcal{I} \subset \mcal{O}_{X}$ be the ideal sheaf of $Z$, and
define
$Z_{n} := \mbf{Spec}\, \mcal{O}_{X} / \mcal{I}^{n+1}$,
$n \geq 0$. The trace maps
$\mcal{K}_{Z_{0}}^{\bdot} \ar \cdots \ar
\mcal{K}_{Z_{n}}^{\bdot} \ar \cdots \ar \mcal{K}_{X}^{\bdot}$
of Proposition \ref{prop1.2} induce a filtration by subcomplexes
$\Gamma_{Z} \mcal{K}_{X}^{\bdot} =
\bigcup_{n=0}^{\infty} \mcal{K}_{Z_{n}}^{\bdot}$. Now
since each morphism $Z_{n} \ar Y$ is proper,
$\opn{Tr}_{f}: f_{*} \mcal{K}_{Z_{n}}^{\bdot} \ar \mcal{K}_{Y}^{\bdot}$
is a homomorphism of complexes.
\end{proof}

\begin{thm} \textup{(Duality)}\ \label{thm2.2}
Let $f: X \ar Y$ be a proper morphism between finite type $k$-schemes.
Then for any complex
$\mcal{M}^{\bdot} \in \mathsf{D}^{\mrm{b}}_{\mrm{c}}(X)$, the
homomorphism
\[ \opn{Hom}_{\mathsf{D}(X)} (\mcal{M}^{\bdot},
\mcal{K}_{X}^{\bdot})
\ar \opn{Hom}_{\mathsf{D}(Y)} (\mrm{R} f_{*} \mcal{M}^{\bdot},
\mcal{K}_{Y}^{\bdot}) \]
induced by
$\opn{Tr}_{f} : f_{*} \mcal{K}_{X}^{\bdot} \ar \mcal{K}_{Y}^{\bdot}$
is an isomorphism.
\end{thm}

\begin{proof}
The proof uses a relative version of Sastry's notion of ``residue
pairs'' and ``pointwise residue pairs'', cf.\ \cite{Ye1} Appendix.
Define a
residue pair relative to $f$ and $\mcal{K}_{Y}^{\bdot}$, to be a pair
$(\mcal{R}^{\bdot}, t)$, with $\mcal{R}^{\bdot}$ a residual
complex on $X$,
and with $t: f_{*} \mcal{R}^{\bdot} \ar \mcal{K}_{Y}^{\bdot}$ a
homomorphism of complexes, which represent the functor
$\mcal{M}^{\bdot} \mapsto \opn{Hom}_{\mathsf{D}(Y)} (\mrm{R} f_{*}
\mcal{M}^{\bdot}, \mcal{K}_{Y}^{\bdot})$ on
$\mathsf{D}^{\mrm{b}}_{\mrm{c}}(X)$.
Such pairs exist; for instance, we may take $\mcal{R}^{\bdot}$ to be the
Cousin complex $f^{\triangle} \mcal{K}_{Y}^{\bdot}$ associated to the
dualizing complex $f^{!} \mcal{K}_{Y}^{\bdot}$
(cf.\ \cite{RD} ch.\ VII \S 3, or ibid.\ Appendix no.\ 4).

A pointwise residue pair relative to $f$ and $\mcal{K}_{Y}^{\bdot}$, is
by definition a pair
$(\mcal{R}^{\bdot}, t)$ as above, but satisfying the condition: for any
closed point $x \in X$, and any coherent $\mcal{O}_{X}$-module $\mcal{M}$
supported on $\{ x \}$, the map
$\opn{Hom}_{\mcal{O}_{X}} (\mcal{M}, \mcal{R}^{\bdot}) \ar
\opn{Hom}_{\mcal{O}_{Y}} (f_{*} \mcal{M}, \mcal{K}_{Y}^{\bdot})$
induced by $t$ is an isomorphism. By the definition of the trace map
$\opn{Tr}_{f}$, the pair
$(\mcal{K}_{X}^{\bdot}, \opn{Tr}_{f})$ is a pointwise residue pair.
In fact,
$k \ar \mcal{O}_{Y, (f(x))} \ar \mcal{O}_{X, (x)}$ are morphisms in
$\mathsf{BCA}(k)$, and by \cite{Ye2} Theorem 7.4 (i),(iv) we get
\[ \opn{Hom}_{\mcal{O}_{X}} (\mcal{M}, \mcal{K}_{X}^{\bdot}) \cong
\opn{Hom}_{\mcal{O}_{Y}} (f_{*} \mcal{M}, \mcal{K}_{Y}^{\bdot}) \cong
\opn{Hom}_{k} (\mcal{M}_{x}, k). \]

The proof of \cite{Ye1} Appendix Theorem 2 goes through also
in the relative situation: the morphism
$\opn{Tr}_{f} : f_{*} \mcal{K}_{X}^{\bdot} \ar \mcal{K}_{Y}^{\bdot}$
in $\mathsf{D}(Y)$ corresponds to a morphism
$\zeta: \mcal{K}_{X}^{\bdot} \ar \mcal{R}^{\bdot}$ in
$\mathsf{D}(X)$. But since both $\mcal{K}_{X}^{\bdot}$ and
$\mcal{R}^{\bdot}$ are residual complexes, $\zeta$ is an actual, unique
homomorphism of complexes (cf.\ \cite{RD} Ch.\ IV Lemma 3.2).
By testing on $\mcal{O}_{X}$-modules
$\mcal{M}$ as above we see that $\zeta$ is indeed an isomorphism of
complexes. So $(\mcal{K}_{X}^{\bdot}, \opn{Tr}_{f})$ is a residue pair.
\end{proof}

Let $\pi: X \ar \opn{Spec} k$ be the structural morphism. In \cite{RD}
\S VII.3 we find the twisted inverse image functor
$\pi^{!} : \mathsf{D}^{+}_{\opn{c}}(k) \ar \mathsf{D}^{+}_{\opn{c}}(X)$.

\begin{cor} \label{cor2.2}
There is a canonical isomorphism
$\zeta_{X} : \mcal{K}_{X}^{\bdot} \iso \pi^{!} k$ in
$\mathsf{D}(X)$. It is compatible with proper and \'{e}tale morphisms.
If $\pi$ is proper then
\[ \opn{Tr}_{\pi} = \opn{Tr}_{\pi}^{\mrm{RD}}
\mrm{R} \pi_{*}(\zeta_{X}) : \pi_{*} \mcal{K}_{X}^{\bdot} \ar k \]
where
$\opn{Tr}_{\pi}^{\mrm{RD}} : \mrm{R} \pi_{*} \pi^{!} k \ar k$
is the trace map of \cite{RD} \S \textup{VII.3}.
\end{cor}

\begin{proof}
The uniqueness of $\zeta_{X}$ follows from considering closed subschemes
$i_{Z} : Z \inj X$ finite over $k$. This is because any endomorphism $a$
of $\mcal{K}_{X}^{\bdot}$ in $\mathsf{D}(X)$
is a global section of $\mcal{O}_{X}$, and $a=1$
iff $i_{Z}^{*}(a) = 1$ for all such $Z$. Existence is proved by covering
$X$ with compactifiable (e.g.\ affine)
open sets and using Theorem \ref{thm2.2}, cf.\
\cite{Ye1} Appendix Theorem 3 and
subsequent Exercise. In particular $\zeta_{X}$ is seen to be compatible
with open immersions. Compatibility with proper morphisms follows from
the transitivity of traces. As for an \'{e}tale morphism $g : U \ar X$,
one has
$g^{*} \mcal{K}_{X}^{\bdot} \cong \mcal{K}_{U}^{\bdot}$ by Theorem
\ref{thm1.1} (b), and also
$g^{*} \pi^{!} k = g^{!} \pi^{!} k = (\pi g)^{!} k$. Testing the
isomorphisms on subschemes $Z \subset U$ finite over $k$ shows that
$g^{*}(\zeta_{X}) = \zeta_{U}$.
\end{proof}

\section{Duals of Differential Operators}

Let $X$ be a $k$-scheme of finite type, where $k$ is a perfect field
of any characteristic.
Suppose $\mcal{M}, \mcal{N}$ are $\mcal{O}_{X}$-modules. By a
differential operator (DO)
$D : \mcal{M} \ar \mcal{N}$ over $\mcal{O}_{X}$,
relative to $k$, we mean
in the sense of \cite{EGA} IV \S 16.8.
Thus $D$ has order $\leq 0$ if $D$ is $\mcal{O}_{X}$-linear,
and
$D$ has order $\leq d$ if for all $a \in \mcal{O}_{X}$,
the commutator $[D, a]$ has order $\leq d-1$.

Recall that the dual of an $\mcal{O}_{X}$-module $\mcal{M}$
is
$\opn{Dual} \mcal{M} =
\mcal{H}om_{\mcal{O}_{X}}^{\bdot}(\mcal{M}, \mcal{K}^{\bdot}_{X})$.
In this section we prove the existence of the dual operator
$\opn{Dual}(D)$, in terms of BCAs and residues. This explicit description
of $\opn{Dual}(D)$ will be needed for the applications in Sections 5-7.
For direct proofs of existence cf.\ Remarks \ref{rem3.2} and
\ref{rem3.3}.

\begin{thm} \label{thm3.1}
Let $\mcal{M},\mcal{N}$ be two $\mcal{O}_{X}$-modules, and let
$D : \mcal{M} \ar \mcal{N}$ be a differential operator of order
$\leq d$. Then there is a homomorphism of graded sheaves
\[ \opn{Dual}(D) :
\opn{Dual} \mcal{N} \ar \opn{Dual} \mcal{M} \]
having the properties below:

\begin{enumerate}
\rmitem{i}  $\opn{Dual}(D)$ is a DO of order $\leq d$.
\rmitem{ii} $\opn{Dual}(D)$ is a homomorphism of complexes.
\rmitem{iii} Functoriality: if
$E : \mcal{N} \ar \mcal{L}$ is another DO, then
$\opn{Dual}(E D) = \opn{Dual}(D)\, \opn{Dual}(E)$.
\rmitem{iv} If $d = 0$, i.e.\ $D$ is $\mcal{O}_{X}$-linear, then
$\opn{Dual}(D)(\phi) = \phi \circ D$ for any
$\phi \in
\mcal{H}om_{\mcal{O}_{X}}^{\bdot}(\mcal{N}, \mcal{K}^{\bdot}_{X})$.
\rmitem{v} Adjunction: under the homomorphisms
$\mcal{M} \ar \opn{Dual} \opn{Dual} \mcal{M}$
and
$\mcal{N} \ar \opn{Dual} \opn{Dual} \mcal{N}$,
one has
$D \mapsto \opn{Dual}(\opn{Dual}(D))$.
\end{enumerate}
\end{thm}

\begin{proof}
By \cite{RD} Theorem II.7.8, an $\mcal{O}_{X}$-module $\mcal{M}'$
is noetherian iff there is a surjection
$\bigoplus_{i = 1}^{n} \mcal{O}_{U_{i}} \surj \mcal{M}'$,
for some open sets $U_{1}, \ldots, U_{n}$. Here $\mcal{O}_{U_{i}}$
is extended by $0$ to a sheaf on $X$.
One has $\mcal{M} = \lim_{\alpha \ar} \mcal{M}_{\alpha}$,
where $\{ \mcal{M}_{\alpha} \}$ is the set of noetherian submodules
of $\mcal{M}$. (We did not assume $\mcal{M},\mcal{N}$ are
quasi-coherent!)
So
\[ \mcal{H}om_{\mcal{O}_{X}}(\mcal{M}, \mcal{K}^{q}_{X})
\cong
\lim_{\leftarrow \alpha}
\mcal{H}om_{\mcal{O}_{X}}(\mcal{M}_{\alpha}, \mcal{K}^{q}_{X}) . \]
Since the sheaf $\mcal{P}^{d}_{X/k}$ of principal parts is coherent,
and
$D : \mcal{M}_{\alpha} \ar \mcal{N}$
induces
\[ \bigoplus_{i = 1}^{n} (\mcal{P}^{d}_{X/k} \otimes \mcal{O}_{U_{i}})
\surj \mcal{P}^{d}_{X/k} \otimes \mcal{M}_{\alpha} \ar \mcal{N} , \]
we conclude that the module
$\mcal{N}_{\alpha} := \mcal{O}_{X} \cdot D(\mcal{M}_{\alpha}) \subset
\mcal{N}$
is also noetherian. Therefore we may assume that both
$\mcal{M}, \mcal{N}$ are noetherian.

We have
\[ \mcal{H}om^{\bdot}_{\mcal{O}_{X}}(\mcal{M}, \mcal{K}^{\bdot}_{X}) =
\bigoplus_{x} \mcal{H}om_{\mcal{O}_{X}}(\mcal{M}, \mcal{K}_{X}(x)) , \]
and
$\mcal{H}om_{\mcal{O}_{X}}(\mcal{M}, \mcal{K}_{X}(x))$
is a constant sheaf with support $\overline{ \{ x \} }$ and module
\[ \mcal{H}om_{\mcal{O}_{X, x}}(\mcal{M}_{x}, \mcal{K}_{X}(x)) =
\opn{Hom}_{A}(M, \mcal{K}(A)) = \opn{Dual}_{A} M , \]
where
$A := \widehat{\mcal{O}}_{X, x} = \mcal{O}_{X, (x)}$
and
$M := A \otimes \mcal{M}_{x}$.
Note that $M$ is a finitely generated $A$-module.
$D : \mcal{M}_{x} \ar \mcal{N}_{x}$
induces a continuous DO
$D : M \ar N = A \otimes \mcal{N}_{x}$
(for the $\mfrak{m}$-adic topology).
According to \cite{Ye2} Theorem 8.6 there is a continuous DO
\[ \opn{Dual}_{A}(D) : \opn{Dual}_{A} N \ar \opn{Dual}_{A} M . \]

Properties (i), (iii), (iv) and (v) follows directly from \cite{Ye2}
Theorem 8.6 and Corollary 8.8.
As for property (ii), consider any saturated chain
$\xi = (x, \ldots, y)$. Since
$\partial^{-} : \mcal{O}_{X, (x)} \ar \mcal{O}_{X, \xi}$ is an
intensification homomorphism, and since
$\partial^{+} : \mcal{O}_{X, (y)} \ar \mcal{O}_{X, \xi}$
is a morphism in $\mathsf{BCA}(k)$ which is also topologically \'{e}tale,
we see that property (ii) is a consequence of \cite{Ye2}
Thm.\ 8.6 and Cor.\ 8.12.
\end{proof}

Let
$\mcal{D}_{X} := \mcal{D}\textit{iff}_{\mcal{O}_{X} / k}(
\mcal{O}_{X}, \mcal{O}_{X})$
be the sheaf of differential operators on $X$. By definition
$\mcal{O}_{X}$ is a left $\mcal{D}_{X}$-module.

\begin{cor} \label{cor3.1}
If $\mcal{M}$ is a left \textup{(}resp.\ right\textup{)}
$\mcal{D}_{X}$-module, then
$\opn{Dual} \mcal{M}$ is a complex of right \textup{(}resp.\
left\textup{)} $\mcal{D}_{X}$-modules. In particular this is true for
$\mcal{K}_{X}^{\bdot} = \opn{Dual} \mcal{O}_{X}$.
\end{cor}

\begin{cor} \label{cor3.2}
Suppose $\mcal{M}^{\bdot}$ is a complex sheaves, where
each $\mcal{M}^{q}$ is an $\mcal{O}_{X}$-module, and
$\mrm{d} : \mcal{M}^{q} \ar \mcal{M}^{q + 1}$ is a DO.
Then there is a dual complex
$\opn{Dual} \mcal{M}^{\bdot}$.
\end{cor}

Specifically,
$(\opn{Dual} \mcal{M}^{\bdot}, \mrm{D})$
is the simple complex associated to the double complex
$(\opn{Dual} \mcal{M}^{\bdot})^{p, q} :=
\mcal{H}om_{\mcal{O}_{X}}(\mcal{M}^{-p}, \mcal{K}_{X}^{q})$.
The operator is
$\mrm{D} = \mrm{D}' + \mrm{D}''$,
where
\[ \begin{split}
\mrm{D}' & := (-1)^{p+q+1} \opn{Dual}(\mrm{d}) :
(\opn{Dual} \mcal{M}^{\bdot})^{p, q} \ar
(\opn{Dual} \mcal{M}^{\bdot})^{p + 1, q} , \\
\mrm{D}'' & := \delta :
(\opn{Dual} \mcal{M}^{\bdot})^{p, q} \ar
(\opn{Dual} \mcal{M}^{\bdot})^{p, q + 1} .
\end{split} \]

It is well known that if $\opn{char} k = 0$ and $X$ is smooth
of dimension $n$, then
$\omega_{X} = \Omega^{n}_{X/k}$ is a right $\mcal{D}_{X}$-module.

\begin{prop} \label{prop3.1}
Suppose $\opn{char} k = 0$ and $X$ is smooth of dimension $n$. Then
$\mrm{C}_{X} : \Omega^{n}_{X/k} \ar \mcal{K}^{-n}_{X}$
\textup{(}the inclusion\textup{)} is $\mcal{D}_{X}$-linear.
\end{prop}

\begin{proof}
It suffice to prove that any $\partial \in \mcal{T}_{X}$
(the tangent sheaf), which we view as a DO
$\partial : \mcal{O}_{X} \ar \mcal{O}_{X}$,
satisfies
$\opn{Dual}(\partial)(\alpha) = - \mrm{L}_{\partial}(\alpha)$,
where $\mrm{L}_{\partial}$
is the Lie derivative, and $\alpha \in \Omega^{n}_{X/k}$.
Localizing at the generic point of $X$ we get
$\partial \in \mcal{D}(k(X))$ and $\alpha \in \omega(k(X))$. Now use
\cite{Ye2} Definition 8.1 and Proposition 4.2.
\end{proof}

\begin{rem} \label{rem3.1}
Proposition \ref{prop3.1} says that in the case
$\opn{char} k = 0$ and $X$ smooth, the $\mcal{D}_{X}$-module
structure on $\mcal{K}^{\bdot}_{X}$ coincides with the standard one,
which is obtained as follows. The quasi-isomorphism
$\mrm{C}_{X} : \Omega^{n}_{X/k} \ar \mcal{K}^{\bdot}_{X}[-n]$
identifies $\mcal{K}^{\bdot}_{X}[-n]$ with the Cousin complex
of $\Omega^{n}_{X/k}$, which is computed in the category
$\msf{Ab}(X)$ (cf.\ \cite{Ha} Section I.2).
Since any $D \in \mcal{D}_{X}$ acts $\mbb{Z}$-linearly on
$\Omega^{n}_{X/k}$, it also acts on $\mcal{K}^{\bdot}_{X}[-n]$.
\end{rem}

\begin{rem} \label{rem3.2}
According to \cite{Sai} there is a direct way to obtain Theorem
\ref{thm3.1} in characteristic $0$. Say $X \subset Y$, with $Y$ smooth.
Then
$\mcal{H}om_{\mcal{O}_{X}}^{\bdot}(\mcal{M}, \mcal{K}^{\bdot}_{X})
\cong
\mcal{H}om_{\mcal{O}_{Y}}^{\bdot}(\mcal{M}, \mcal{K}^{\bdot}_{Y})$.
Now by \cite{Sai} \S 2.2.3
any DO $D : \mcal{M} \ar \mcal{N}$ of order $\leq d$ can be
viewed as
\[ D \in
\mcal{H}om_{\mcal{O}_{Y}}(\mcal{M}, \mcal{N} \otimes_{\mcal{O}_{Y}}
\mcal{D}^{d}_{Y}) \subset
\mcal{H}om_{\mcal{D}_{Y}}(
\mcal{M} \otimes_{\mcal{O}_{Y}} \mcal{D}_{Y},
\mcal{N} \otimes_{\mcal{O}_{Y}} \mcal{D}_{Y})  \]
(right $\mcal{D}_{Y}$-modules). Since $\mcal{K}^{q}_{Y}$ is a
$\mcal{D}_{Y}$-module (cf.\ Remark \ref{rem3.1}), we get
\[ \mcal{H}om_{\mcal{O}_{Y}}(\mcal{M}, \mcal{K}^{q}_{Y})
\cong
\mcal{H}om_{\mcal{D}_{Y}}(\mcal{M} \otimes_{\mcal{O}_{Y}}
\mcal{D}_{Y}, \mcal{K}^{q}_{Y}) \]
and so we obtain the dual operator $\opn{Dual}(D)$.
I thank the referee for pointing out this fact to me.
\end{rem}

\begin{rem} \label{rem3.3}
Suppose $\opn{char} k = p > 0$. Then a $k$-linear map
$D : \mcal{M} \ar \mcal{N}$ is a DO over $\mcal{O}_{X}$ iff it is
$\mcal{O}_{X^{(p^{n}/k)}}$-linear, for $n \gg 0$.
Here
$X^{(p / k)} \ar X$ is the Frobenius morphism relative to $k$,
cf.\ \cite{Ye1} Theorem 1.4.9.
Since $\opn{Tr}$ induces an isomorphism
\[ \mcal{H}om_{\mcal{O}_{X}}^{\bdot}(\mcal{M}, \mcal{K}^{\bdot}_{X})
\cong
\mcal{H}om_{\mcal{O}_{X^{(p^{n}/k)}}}^{\bdot}(
\mcal{M}, \mcal{K}^{\bdot}_{X^{(p^{n}/k)}}) \]
we obtain the dual operator $\opn{Dual}(D)$.
\end{rem}

Let us finish this section with an application to rings of differential
operators. Given a finitely generated (commutative) $k$-algebra $A$,
denote by $\mcal{D}(A) := \opn{Diff}_{A/k}(A,A)$
the ring of differential operators over $A$. Such rings are of interest
for ring theorists (cf.\ \cite{MR} and \cite{HoSt}).
It is well known that if $\opn{char} k = 0$ and $A$ is smooth,
then the opposite ring
$\mcal{D}(A)^{\circ} \cong \omega_{A} \otimes_{A} \mcal{D}(A)
\otimes_{A} \omega_{A}^{-1}$,
where $\omega_{A} = \Omega^{n}_{A / k}$.
The next theorem is a vast generalization of this fact.

Given complexes $M^{\bdot}, N^{\bdot}$ of $A$-modules let
$\opn{Diff}^{\bdot}_{A / k}(M^{\bdot}, N^{\bdot})$
be the complex of $k$-modules  which in degree $n$ is
$\prod_{p} \opn{Diff}_{A / k}(M^{p}, N^{p+n})$.
Let
$\mcal{K}^{\bdot}_{A} := \Gamma(X, \mcal{K}^{\bdot}_{X})$
with
$X := \opn{Spec} A$.
By Corollary \ref{cor3.1}, it is a complex of right $\mcal{D}(A)$-modules.

\begin{thm} \label{thm3.2}
There is a natural isomorphism of filtered $k$-algebras
\[ \mcal{D}(A)^{\circ} \cong
\opn{H}^{0} \opn{Diff}_{A / k}^{\bdot}(
\mcal{K}_{A}^{\bdot}, \mcal{K}_{A}^{\bdot}) . \]
\end{thm}

\begin{proof}
First observe that since DOs preserve support,
$\opn{Diff}_{A / k}(\mcal{K}_{A}^{p}, \mcal{K}_{A}^{p-1}) = 0$
for all $p$. This means that every local section
$D \in \opn{H}^{0} \opn{Diff}_{A / k}^{\bdot}(\mcal{K}_{A}^{\bdot},
\mcal{K}_{A}^{\bdot})$
is a well defined DO
$D : \mcal{K}_{A}^{\bdot} \ar \mcal{K}_{A}^{\bdot}$
which commutes with the coboundary $\delta$. Applying
$\opn{Dual}$ and taking $0$-th cohomology we obtain a DO
\[ D^{\vee} = \opn{H}^{0} \opn{Dual}(D) :
\opn{H}^{0} \opn{Dual} \mcal{K}_{A}^{\bdot} \ar
\opn{H}^{0} \opn{Dual} \mcal{K}_{A}^{\bdot} . \]
But
$\opn{H}^{0} \opn{Dual} \mcal{K}_{A}^{\bdot} = A$,
so $D^{\vee} \in \mcal{D}(A)$.
Finally according to Theorem \ref{thm3.1} (v), $D = D^{\vee \vee}$
for
$D \in \opn{H}^{0} \opn{Diff}_{A / k}^{\bdot}(\mcal{K}_{A}^{\bdot},
\mcal{K}_{A}^{\bdot})$
or $D \in \mcal{D}(A)$.
\end{proof}

Recall that an $n$-dimensional integral domain $A$ is a Gorenstein algebra
iff
$\omega_{A} = \mrm{H}^{-n} \mcal{K}_{A}^{\bdot} \ar
\mcal{K}_{A}^{\bdot}[-n]$
is a quasi-isomorphism, and $\omega_{A}$ is invertible.

\begin{cor} \label{cor3.6}
If $A$ is a Gorenstein $k$-algebra,
there is a canonical isomorphism of filtered $k$-algebras
\[ \mcal{D}(A)^{\circ} \cong
\opn{Diff}_{A / k}(\omega_{A}, \omega_{A}) \cong
\omega_{A} \otimes_{A} \mcal{D}(A) \otimes_{A} \omega_{A}^{-1} . \]
\end{cor}

\begin{rem}
In \cite{Ho}, the right $\mcal{D}(A)$-module structure on $\omega_{A}$
was exhibited, when $X = \opn{Spec} A$ is a curve. Corollary
\ref{cor3.6} was proved there for complete intersection curves.
\end{rem}

\section{The De Rham-Residue Complex}

As before $k$ is a perfect field of any characteristic.
Let $X$ be a $k$-scheme of finite type.
In this section we define a canonical complex on $X$, the De
Rham-residue complex $\mcal{F}^{\bdot}_{X}$.
As we shall see in Corollary \ref{cor4.3}, $\mcal{F}^{\bdot}_{X}$
coincides (up to indices and signs) with the double complex
$\mcal{K}^{\bdot, *}_{X}$ of \cite{EZ}.

According to Theorem \ref{thm3.1}, if $\mcal{M}^{\bdot}$ is a complex
of sheaves, with each $\mcal{M}^{q}$ an $\mcal{O}_{X}$-module and
$\mrm{d} : \mcal{M}^{q} \ar \mcal{M}^{q+1}$
a DO, then $\opn{Dual} \mcal{M}^{\bdot}$ is a complex of the same kind.

\begin{dfn} \label{dfn4.1}
The {\em De Rham-residue complex} on $X$ is the complex
\[ \mcal{F}_{X}^{\bdot} :=
\opn{Dual} \Omega^{\bdot}_{X / k} . \]
of Corollary \ref{cor3.2}.
\end{dfn}

Note that the double complex $\mcal{F}_{X}^{\bdot \bdot}$ is
concentrated in the third quadrant of the $(p,q)$-plane.

\begin{prop} \label{prop4.2}
$\mcal{F}_{X}^{\bdot}$ is a right DG module over
$\Omega^{\bdot}_{X/k}$.
\end{prop}

\begin{proof}
The graded module structure is clear.
It remains to check that
\[ \mrm{D}(\phi \cdot \alpha) = (\mrm{D} \phi) \cdot \alpha +
(-1)^{p+q} \phi \cdot (\mrm{d} \alpha) \]
for $\phi \in \mcal{F}_{X}^{p,q}$ and
$\alpha \in \Omega^{p'}_{X/k}$.
But this is a straightforward computation using Theorem \ref{thm3.1}.
\end{proof}

\begin{prop} \label{prop4.3}
Let $g: U \ar X$ be \'{e}tale. Then there is a homomorphism of
complexes
$\mrm{q}_{g} : \mcal{F}_{X}^{\bdot} \ar g_{*} \mcal{F}_{U}^{\bdot}$,
which induces an isomorphism of graded sheaves
$1 \otimes \mrm{q}_{g} : g^{*} \mcal{F}_{X}^{\bdot} \iso
\mcal{F}_{U}^{\bdot}$.
\end{prop}

\begin{proof}
Consider the isomorphisms
$g^{*} \Omega^{\bdot}_{X/k} \cong \Omega^{\bdot}_{U/k}$
and
$1 \otimes \mrm{q}_{g} : g^{*} \mcal{K}_{X}^{\bdot} \ar
\mcal{K}_{U}^{\bdot}$
of Theorem \ref{thm1.1}. Clearly
$1 \otimes \mrm{q}_{g} :
g^{*} \mcal{F}_{X}^{p,q} \ar \mcal{F}_{U}^{p,q}$
is an isomorphism.
In light of \cite{Ye2} Theorem 8.6 (iv),
$\mrm{q}_{g} : \mcal{F}_{X}^{\bdot} \ar g_{*} \mcal{F}_{U}^{\bdot}$
is a homomorphism of complexes.
\end{proof}

Let $f: X \ar Y$ be a morphism of schemes. Define a homomorphism of
graded sheaves
$\opn{Tr}_{f} : f_{*} \mcal{F}_{X}^{\bdot} \ar
\mcal{F}_{Y}^{\bdot}$
by composing
$f^{*} : \Omega^{\bdot}_{Y/k} \ar f_{*} \Omega^{\bdot}_{X/k}$
with
$\opn{Tr}_{f} : f_{*} \mcal{K}_{X}^{\bdot} \ar
\mcal{K}_{Y}^{\bdot}$
of Definition \ref{dfn1.3}.

\begin{prop} \label{prop4.1}
$\opn{Tr}_{f}$ commutes with $\mrm{D}'$. If $f$ is proper then
$\opn{Tr}_{f}$ also commutes with $\mrm{D}''$.
\end{prop}

\begin{proof}
Let $y \in Y$ and let $x$ be a closed point in $f^{-1}(y)$. Then
$f^{*}: \mcal{O}_{Y,(y)} \ar \mcal{O}_{X,(x)}$ is a morphism in
$\msf{BCA}(k)$. Applying \cite{Ye2} Cor.\ 8.12 to the DOs
\[ \mrm{d} f^{*} = f^{*} \mrm{d} :
\Omega^{p}_{Y/k, (y)} \ar \Omega^{p+1}_{X/k, (x)} \]
we get a dual homomorphism
\[ \opn{Dual}_{f^{*}}(\mrm{d} f^{*}) =
\opn{Dual}_{f^{*}}(f^{*} \mrm{d}) : \opn{Dual}_{\mcal{O}_{X,(x)}}
\Omega^{p+1}_{X/k, (x)} \ar
\opn{Dual}_{\mcal{O}_{Y,(y)}} \Omega^{p}_{Y/k, (y)} , \]
which equals both
$\opn{Tr}_{f} \opn{Dual}_{X}(\mrm{d})$ and
$\opn{Dual}_{Y}(\mrm{d}) \opn{Tr}_{f}$.
The commutation of $\mrm{D}''$ with $\opn{Tr}_{f}$ in the proper case
is immediate from Thm.\ \ref{thm2.1}.
\end{proof}

Of course if $f : X \ar Y$ is a closed immersion, then $\opn{Tr}_{f}$ is
injective, and it identifies $\mcal{F}_{X}^{\bdot}$ with the
subsheaf
$\mcal{H}om_{\Omega^{\bdot}_{Y / k}}(\Omega^{\bdot}_{X / k},
\mcal{F}_{Y}^{\bdot})$ of
$\mcal{F}_{Y}^{\bdot}$.
Just as in Corollary \ref{cor2.1} we get:

\begin{cor} \label{cor4.1}
Let $f : X \ar Y$ be a morphism of schemes, and
let $Z \subset X$ be a closed subscheme which is proper over $Y$. Then
the trace map
$\opn{Tr}_{f}: f_{*} \Gamma_{Z} \mcal{F}_{X}^{\bdot} \ar
\mcal{F}_{Y}^{\bdot}$
is a homomorphism of complexes.
\end{cor}

Suppose $X$ is an integral scheme of dimension $n$. The canonical
homomorphism
\begin{equation}
\opn{C}_{X} : \Omega^{n}_{X/k} \ar
\mcal{K}_{X}^{-n} = k(X) \otimes_{\mcal{O}_{X}} \Omega^{n}_{X/k}
\end{equation}
can be viewed as a global section of $\mcal{F}_{X}^{-n, -n}$.

\begin{lem} \label{lem4.3}
Suppose $X$ is an integral scheme. Then
$\mrm{D}' (\opn{C}_{X}) = \mrm{D}'' (\opn{C}_{X}) = 0$.
\end{lem}

\begin{proof}
By \cite{Ye1} Section 4.5,
$\mrm{D}''(\opn{C}_{X}) = \pm \delta (\opn{C}_{X}) = 0$.
Next, let $K := k(X)$. Choose $t_{1}, \ldots t_{n} \in K$
such that
$\Omega^{1}_{K / k} = \bigoplus K \cdot \mrm{d} t_{i}$.
Taking products of the $\mrm{d} t_{i}$ as bases of
$\Omega^{n - 1}_{K / k}$ and $\Omega^{n}_{K / k}$, we see from
\cite{Ye2} Theorem 8.6 and Definition 8.1 that
$\opn{Dual}_{K}(\opn{C}_{X}) = 0$.
\end{proof}

\begin{prop} \label{prop4.4}
If $X$ is smooth irreducible of dimension $n$, then the DG homomorphism
$\Omega^{\bdot}_{X/k} \ar \mcal{F}^{\bdot}_{X}[-2n]$,
$\alpha \mapsto \mrm{C}_{X} \cdot \alpha$,
is a quasi-isomorphism.
\end{prop}

\begin{proof}
First note that $\mrm{D}(\mrm{C}_{X}) = 0$, so this is indeed a
DG homomorphism.
Filtering these complexes according to the $p$-degree
we reduce to looking at
$\Omega^{p}_{X/k}[n] \ar \mcal{F}^{p-n, \bdot}_{X}$.
That is a quasi-isomorphism by Theorem \ref{thm1.1} part d.
\end{proof}

\begin{cor} \label{cor4.3}
The complex $\mcal{F}^{\bdot}_{X}$ is the same as the complex
$\mcal{K}^{\bdot, *}_{X}$ of \cite{EZ}, up to signs and indexing.
\end{cor}

\begin{proof}
If $X$ is smooth of dimension $n$ this is because
$\mcal{F}^{\bdot}_{X} \cong
\Omega^{\bdot}_{X/k}[n] \otimes \mcal{K}^{\bdot}_{X}$
is the Cousin complex of
$\bigoplus \Omega^{p}_{X/k}[p]$, and $\mrm{D}'$
is (up to sign) the Cousin functor applied to $\mrm{d}$.
If $X$ is a general scheme embedded in a smooth scheme $Y$,
use Proposition \ref{prop4.1}.
\end{proof}

\begin{dfn} \label{dfn4.4}
Given a scheme $X$, let $X_{1}, \ldots, X_{r}$ be its irreducible
components, with their induced reduced subscheme structures. For each $i$
let $x_{i}$ be the generic point of $X_{i}$, and let $f_{i} : X_{i} \ar X$
be the inclusion morphism. We define the fundamental class $\opn{C}_{X}$
by:
\[ \opn{C}_{X}:= \sum_{i=1}^{r} \opn{length}(\mcal{O}_{X, x_{i}})
\opn{Tr}_{f_{i}}(\opn{C}_{X_{i}}) \in \mcal{F}_{X}^{\bdot} . \]
\end{dfn}

The next proposition is easily verified using Propositions \ref{prop4.1}
and \ref{prop4.3}. It should be compared to
\cite{EZ}  Theorem III.3.1.

\begin{prop} \label{prop4.7}
For any scheme $X$, the fundamental class
$\opn{C}_{X} \in$ \linebreak
$\Gamma(X, \mcal{F}_{X}^{\bdot})$
is annihilated
by $\mrm{D}'$ and $\mrm{D}''$. If $X$ has pure dimension $n$,
then $\opn{C}_{X}$ has bidegree $(-n,-n)$. If $f: X \ar Y$ is a proper,
surjective, generically finite morphism between integral schemes, then
$\opn{Tr}_{f}(\opn{C}_{X}) = \opn{deg}(f) \opn{C}_{Y}$.
If $g: U \ar X$ is \'{e}tale, then
$\mrm{C}_{U} = \mrm{q}_{g}(\mrm{C}_{X})$.
\end{prop}

\begin{rem}
In \cite{Ye4} it is shown that $\mcal{F}_{X}^{\bdot}$ is a right
DG module over the DGA of Beilinson adeles
$\mcal{A}_{X}^{\bdot} = \ul{\mbb{A}}^{\bdot}_{\mrm{red}}(
\Omega^{\bdot}_{X / k})$.
Now let $\mcal{E}$ be a locally free $\mcal{O}_{X}$-module of rank
$r$, and let $Z \subset X$ be the zero locus of a regular section of
$\mcal{E}$. According to the adelic Chern-Weil theory of \cite{HY2}
there is an adelic connection $\nabla$ on $\mcal{E}$ such that
the Chern form
$\mrm{c}_{r}(\mcal{E}; \nabla) \in \mcal{A}_{X}^{2r}$
satisfies
$\mrm{C}_{Z} = \pm \mrm{C}_{X} \cdot \mrm{c}_{r}(\mcal{E}; \nabla)
\in \mcal{F}_{X}^{\bdot}$.
\end{rem}

\section{De Rham Homology and the Niveau Spectral Sequence}

Let $X$ be a finite type scheme over a field $k$ of characteristic $0$.
In \cite{Ye3} it is shown that if $X \subset \mfrak{X}$ is
a smooth formal embedding (see below) then the De Rham complex
$\widehat{\Omega}^{\bdot}_{\mfrak{X} / k}$ calculates the De Rham
cohomology $\mrm{H}^{\bdot}_{\mrm{DR}}(X)$.
In this section we will show that the De Rham-residue complex
$\mcal{F}^{\bdot}_{\mfrak{X}}$ of $\mfrak{X}$ calculates
the De Rham homology $\mrm{H}^{\mrm{DR}}_{\bdot}(X)$.
This is done by computing the niveau spectral sequence
converging to $\mrm{H}^{\bdot}(X, \mcal{F}^{\bdot}_{\mfrak{X}})$
(Theorem \ref{thm5.3}).
We will draw a few conclusions, including
the contravariance of homology w.r.t.\ \'{e}tale morphisms
(Theorem \ref{thm5.1}).
As a reference for algebraic De Rham (co)homology we suggest
\cite{Ha}.

Given a noetherian adic formal scheme $\mfrak{X}$ and a defining ideal
$\mcal{I} \subset \mcal{O}_{\mfrak{X}}$,
let $X_{n}$ be the (usual) noetherian scheme
$(\mfrak{X}, \mcal{O}_{\mfrak{X}} / \mcal{I}^{n+1})$.
Suppose $f : \mfrak{X} \ar \mfrak{Y}$ is a morphism between such
formal schemes, and let
$\mcal{I} \subset \mcal{O}_{\mfrak{X}}$ and
$\mcal{J} \subset \mcal{O}_{\mfrak{Y}}$
be defining ideals such that
$f^{-1} \mcal{J} \cdot \mcal{O}_{\mfrak{X}} \subset \mcal{I}$.
Such ideals are always available.
We get a morphism of (usual) schemes
$f_{0} : X_{0} \ar Y_{0}$.

\begin{dfn} \label{dfn5.2}
A morphism $f : \mfrak{X} \ar \mfrak{Y}$ between (noetherian)
adic formal schemes is called {\em formally finite type}
(resp.\ {\em formally finite} or {\em formally proper})
if the morphism $f_{0} : X_{0} \ar Y_{0}$
is finite type (resp.\ finite or proper).
\end{dfn}

Obviously these notions are independent of the particular defining
ideals chosen.

\begin{exa} \label{exa5.1}
If $X \ar Y$ is a finite type morphism of noetherian schemes,
$X_{0} \subset X$ is a locally closed subset and
$\mfrak{X} = X_{/ X_{0}}$ is the formal completion, then
$\mfrak{X} \ar Y$ is formally finite type.
Such a morphism is called {\em algebraizable}.
\end{exa}

\begin{dfn}
A morphism of formal schemes $\mfrak{X} \ar \mfrak{Y}$ is said to be {\em
formally smooth} (resp.\ {\em formally \'{e}tale}) if, given
a (usual) affine scheme $Z$ and a closed subscheme $Z_{0} \subset Z$
defined by a nilpotent ideal, the map
$\opn{Hom}_{\mfrak{Y}}(Z, \mfrak{X})$ \linebreak
$\ar \opn{Hom}_{\mfrak{Y}}(Z_{0}, \mfrak{X})$
is surjective (resp.\ bijective).
\end{dfn}

This is the definition of formal smoothness used in \cite{EGA} IV
Section 17.1. We shall also require the next notion.

\begin{dfn}
A morphism $g: \mfrak{X} \ar \mfrak{Y}$ between noetherian formal
sche\-mes is called {\em \'{e}tale} if it is of finite
type (see \cite{EGA} I \S 10.13) and formally \'{e}tale.
\end{dfn}

Note that if $\mfrak{Y}$ is a usual scheme, then so is
$\mfrak{X}$, and $g$ is an \'{e}tale morphism of schemes.

\begin{dfn} \label{dfn5.3}
A {\em smooth formal embedding} of $X$ (over $k$)
is a closed immersion of $X$ into a formal scheme $\mfrak{X}$,
which induces a homeomorphism on the underlying topological spaces, and
such that
$\mfrak{X}$ is of formally finite type and formally smooth over $k$.
\end{dfn}

\begin{exa}
If $X$ is smooth over $Y = \opn{Spec} k$ and $X_{0}, \mfrak{X}$
are as in Example \ref{exa5.1}, then $X_{0} \subset \mfrak{X}$
is a smooth formal embedding.
\end{exa}

Let $\xi = (x_{0}, \ldots, x_{q})$ be a saturated chain of points
in the formal scheme $\mfrak{X}$. Choose a defining ideal
$\mcal{I}$, and let $X_{n}$ be as above.
Define the Beilinson completion
$\mcal{O}_{\mfrak{X}, \xi} := \lim_{\leftarrow n}
\mcal{O}_{X_{n}, \xi}$
(which of course is independent of $\mcal{I}$).

\begin{lem} \label{lem5.2}
Let $\mfrak{X}$ be formally finite type over $k$, and let
$\xi$ be a saturated chain in $\mfrak{X}$.
Then $\mcal{O}_{\mfrak{X}, \xi}$ is a BCA over $k$. If
$\mfrak{X} = X_{/ X_{0}}$, then
$\mcal{O}_{\mfrak{X}, \xi} \cong \mcal{O}_{X, \xi}$.
\end{lem}

\begin{proof}
First assume $\mfrak{X} = X_{/ X_{0}}$. Taking $\mcal{I}$ to be
the ideal of $X_{0}$ in $X$, we have
\begin{multline*}
\hspace{1cm}
\mcal{O}_{\mfrak{X}, \xi} =
\lim_{\leftarrow n} (\mcal{O}_{X} / \mcal{I}^{n})_{\xi} \cong
\lim_{\leftarrow m,n} \mcal{O}_{X, \xi} /
(\mcal{I}^{n} \mcal{O}_{X, \xi} + \mfrak{m}_{\xi}^{m}) \\
\cong \lim_{\leftarrow m} \mcal{O}_{X, \xi} / \mfrak{m}_{\xi}^{m} =
\mcal{O}_{X, \xi} . \hspace{1cm}
\end{multline*}
Now by \cite{Ye3} Proposition 1.20 and Lemma 1.1, locally there is a
closed immersion $\mfrak{X} \subset \mfrak{Y}$, with
$\mfrak{Y}$ algebraizable (i.e.\ $\mfrak{Y} = Y_{/ Y_{0}}$).
So there is a surjection
$\mcal{O}_{\mfrak{Y}, \xi} \ar
\mcal{O}_{\mfrak{X}, \xi}$,
and this implies that $\mcal{O}_{\mfrak{X}, \xi}$ is a BCA.
\end{proof}

One can construct the complexes
$\mcal{K}_{\mfrak{X}}^{\bdot}$ and
$\mcal{F}_{\mfrak{X}}^{\bdot}$ for a formally finite type
formal scheme $\mfrak{X}$, as follows. Define
$\mcal{K}_{\mfrak{X}}(x) :=
\mcal{K}(\mcal{O}_{\mfrak{X}, (x)})$.
Now let $(x,y)$ be a saturated chain.
Then there is an intensification homomorphism
$\partial^{-}: \mcal{O}_{\mfrak{X}, (x)} \ar
\mcal{O}_{\mfrak{X}, (x,y)}$
and a morphism of BCAs
$\partial^{+}: \mcal{O}_{\mfrak{X}, (y)} \ar
\mcal{O}_{\mfrak{X}, (x,y)}$.
Therefore we get a homomorphism of
$\mcal{O}_{\mfrak{X}}$-modules
$\delta_{(x,y)}: \mcal{K}_{\mfrak{X}}(x) \ar
\mcal{K}_{\mfrak{X}}(y)$. Define a graded sheaf
$\mcal{K}_{\mfrak{X}}^{\bdot} =
\bigoplus_{x \in \mfrak{X}} \mcal{K}_{\mfrak{X}}(x)$
on $\mfrak{X}$, as in \S 1.
Let
$\widehat{\Omega}^{\bdot}_{\mfrak{X}/k}$  be the complete De Rham
complex on $\mfrak{X}$, and set
$\mcal{F}_{\mfrak{X}}^{p, q} :=
\mcal{H}om_{\mcal{O}_{\mfrak{X}}}(
\widehat{\Omega}^{-p}_{\mfrak{X}/k},
\mcal{K}_{\mfrak{X}}^{q})$.

\begin{prop} \label{prop5.5}
Let $\mfrak{X}$ be a formally finite type formal scheme over $k$.
\begin{enumerate}
\item $\mcal{F}_{\mfrak{X}}^{\bdot}$
is a complex.
\item If $g : \mfrak{U} \ar \mfrak{X}$ is \'{e}tale, then
there is a homomorphism of complexes
$\mrm{q}_{g} : \mcal{F}_{\mfrak{X}}^{\bdot} \ar g_{*}
\mcal{F}_{\mfrak{U}}^{\bdot}$,
which induces an isomorphism of graded sheaves
$1 \otimes \mrm{q}_{g} : g^{*} \mcal{F}_{\mfrak{X}}^{\bdot} \iso
\mcal{F}_{\mfrak{U}}^{\bdot}$.
\item If $f : \mfrak{X} \ar \mfrak{Y}$ is formally proper, then
there is a homomorphism of complexes
$\opn{Tr}_{f} : f^{*} \mcal{F}_{\mfrak{X}}^{\bdot} \ar
\mcal{F}_{\mfrak{Y}}^{\bdot}$.
\end{enumerate}
\end{prop}

\begin{proof}
1.\ Let $X_{n} \subset \mfrak{X}$ be as before. Then one has
$\mcal{F}_{\mfrak{X}}^{\bdot} =
\bigcup \mcal{F}_{X_{n}}^{\bdot}$, so this is a complex.\\
2.\ Take $U_{n} := \mfrak{U} \times_{\mfrak{X}}
X_{n}$; then each $U_{n} \ar X_{n}$ is an \'{e}tale morphism of schemes,
and we can use Proposition \ref{prop4.3}.\\
3.\ Apply Proposition \ref{prop4.1} to $X_{n} \ar Y_{n}$.
\end{proof}

\begin{prop} \label{prop5.1}
Assume $\mfrak{X} = Y_{/X}$ for some smooth irreducible scheme
$Y$ of dimension $n$ and closed set $X \subset Y$. Then there is
a natural isomorphism of complexes
\begin{equation} \label{eqn5.3}
\mcal{F}^{\bdot}_{\mfrak{X}} \cong
\ul{\Gamma}_{X} \mcal{F}^{\bdot}_{Y} .
\end{equation}
Hence
$\mcal{F}^{\bdot}_{\mfrak{X}} \cong
\mrm{R} \ul{\Gamma}_{X} \Omega^{\bdot}_{Y / k}[2n]$
in the derived category $\mathsf{D}(\mathsf{Ab}(Y))$, and consequently
\[ \mrm{H}^{-q}(X, \mcal{F}^{\bdot}_{\mfrak{X}}) \cong
\mrm{H}^{2n-q}_{X}(Y, \Omega^{\bdot}_{Y / k}) =
\mrm{H}^{\mrm{DR}}_{q}(X) . \]
\end{prop}

\begin{proof}
The isomorphism (\ref{eqn5.3}) is immediate from Lemma \ref{lem5.2}.
But according to Proposition \ref{prop4.4}, $\mcal{F}^{\bdot}_{Y}$
is a flasque resolution of $\Omega^{\bdot}_{Y/k}[2n]$
in $\mathsf{Ab}(Y)$.
\end{proof}

We need some algebraic results, phrased in the terminology of
\cite{Ye1} \S 1.
Let $K$ be a complete, separated semi-topological (ST) commutative
$k$-algebra, and let
$\ul{t} = (t_{1}, \ldots, t_{n})$ be a sequence of indeterminates.
Let $K[[\, \ul{t}\, ]]$ and $K((\ul{t}))$ be the rings of formal
power series, and of iterated Laurent series, respectively. These
are complete, separated ST $k$-algebras.
Let $T$ be the free $k$-module with basis
$\alpha_{1}, \ldots, \alpha_{n}$ and let
$\bigwedge_{k}^{\bdot} T$ be the exterior algebra over $k$.

\begin{lem} \label{lem5.1} \textup{(``Poincar\'{e} Lemma'')}\
The DGA homomorphisms
\[ \Omega^{\bdot, \opn{sep}}_{K / k} \ar
\Omega^{\bdot, \opn{sep}}_{K[[\, \ul{t}\, ]] / k} \]
and
\[ \Omega^{\bdot, \opn{sep}}_{K / k} \otimes_{k}
\bigwedge\nolimits_{k}^{\bdot} T
\ar
\Omega^{\bdot, \opn{sep}}_{K((\ul{t})) / k},\
\alpha_{i} \mapsto \opn{dlog} t_{i} \]
are quasi-isomorphisms.
\end{lem}

\begin{proof}
Since
\[ \Omega^{\bdot, \opn{sep}}_{K[[\, \ul{t}\, ]] / k} \cong
K[[\, \ul{t}\, ]] \otimes_{k[\, \ul{t}\, ]}
\Omega^{\bdot}_{k[\, \ul{t}\, ] / k} \]
the homotopy operator (``integration'') of the Poincar\'{e} Lemma
for the graded polynomial algebra $k[\, \ul{t}\, ]$ works here also.

For $K((t))$ (i.e.\ $n = 1$) we have
\[ \Omega^{\bdot, \opn{sep}}_{K((t)) / k} \cong
\Omega^{\bdot, \opn{sep}}_{K[[t]] / k} \oplus
\Omega^{\bdot, \opn{sep}}_{K[t^{-1}] / k} \wedge \opn{dlog} t \]
so we have a quasi-isomorphism. For $n > 1$ use
induction on $n$ and the  K\"{u}nneth formula.
\end{proof}

\begin{lem} \label{lem5.4}
Suppose $A$ is a local BCA and
$\sigma, \sigma' : K \ar A$ are two coefficient fields. Then
\[ \mrm{H}(\sigma) = \mrm{H}(\sigma') :
\mrm{H} \Omega^{\bdot, \opn{sep}}_{K/k} \ar
\mrm{H} \Omega^{\bdot, \opn{sep}}_{A/k} . \]
\end{lem}

\begin{proof}
Choosing generators for the maximal ideal of $A$, $\sigma$
induces a surjection of BCAs $\tilde{A} = K [[\, \ul{t}\, ]] \ar A$.
Denote by $\tilde{\sigma} : K \ar \tilde{A}$ the inclusion.
The coefficient field $\sigma'$ lifts to some coefficient field
$\tilde{\sigma}' : K \ar \tilde{A}$. It suffices to show that
$\mrm{H}(\tilde{\sigma}) = \mrm{H}(\tilde{\sigma}')$.
But by Lemma \ref{lem5.1} both of these are bijective, and using the
projection $\tilde{A} \ar K$ we see they are in fact equal.
\end{proof}

Given a saturated chain $\xi = (x, \ldots, y)$ in $X$
and a coefficient field
$\sigma: k(y) \ar \mcal{O}_{X, (y)}$, there is the {\em Parshin residue
map}
\[ \opn{Res}_{\xi, \sigma}: \Omega^{\bdot}_{k(x)/k} \ar
\Omega^{\bdot}_{k(y)/k} \]
(cf.\ \cite{Ye1} Definition 4.1.3). It is a
map of DG $k$-modules of degree equal to $-(\text{length of } \xi)$.

\begin{prop} \label{prop5.2}
Let $\xi = (x, \ldots, y)$ be a saturated chain in $X$. Then
the map of graded $k$-modules
\[ \opn{Res}_{\xi} := \mrm{H}(\opn{Res}_{\xi, \sigma}) :
\mrm{H} \Omega^{\bdot}_{k(x)/k} \ar
\mrm{H} \Omega^{\bdot}_{k(y)/k} \]
is independent of the coefficient field $\sigma$.
\end{prop}

\begin{proof}
Say $\xi$ has length $n$.
Let $L$ be one of the local factors of $k(\xi) = k(x)_{\xi}$,
so it is an $n$-dimensional topological local field (TLF).
Let $K := \kappa_{n}(L)$, the last residue field of $L$,
which is a finite separable $k(y)$-algebra.
Then $\sigma$ extends uniquely to a morphism of TLFs
$\sigma: K \ar L$, and it is certainly enough to check that
\begin{equation} \label{eqn5.4}
\mrm{H}(\opn{Res}_{L/K; \sigma}) :
\mrm{H} \Omega^{\bdot, \mrm{sep}}_{L/k} \ar
\mrm{H} \Omega^{\bdot, \mrm{sep}}_{K/k}
\end{equation}
is independent of $\sigma$.

After choosing a regular system of parameters
$\ul{t} = (t_{1}, \ldots, t_{n})$ in $L$ we get
$L \cong K((\ul{t}))$. According to Lemma \ref{lem5.1},
$\mrm{H}(\sigma)$ induces an isomorphism of $k$-algebras
\begin{equation} \label{eqn5.5}
\mrm{H} \Omega^{\bdot}_{K/k} \otimes_{k}
\bigwedge\nolimits_{k}^{\bdot} T \cong
\mrm{H} \Omega^{\bdot, \mrm{sep}}_{L / k} .
\end{equation}
But by Lemma \ref{lem5.4} this isomorphism is independent of $\sigma$.
The map (\ref{eqn5.4}) is
$\mrm{H} \Omega^{\bdot}_{K/k}$-linear, and
it sends $\bigwedge^{p}_{k} T$ to $0$ if $p < n$, and
$\opn{dlog} t_{1} \wedge \cdots \wedge \opn{dlog} t_{n}
\mapsto 1$. Hence (\ref{eqn5.4}) is independent of $\sigma$.
\end{proof}

The topological space $X$ has an increasing filtration by families of
supports
$\emptyset = X_{-1} \subset X_{0} \subset X_{1} \subset
\cdots$,
where
\[ X_{q} := \{ Z \subset X \mid Z \text{ is closed and }
\opn{dim} Z \leq q \} . \]
We write
$x \in X_{q} / X_{q-1}$ if $\overline{\{x\}} \in X_{q} - X_{q-1}$,
and the set $X_{q} / X_{q-1}$ is called the $q$-skeleton of $X$.
(This notation is in accordance with \cite{BlO}; in \cite{Ye1} $X_{q}$
denotes the $q$-skeleton.)
The {\em niveau filtration} on $\mcal{F}_{\mfrak{X}}^{\bdot}$ is
$\mrm{N}_{q} \mcal{F}_{\mfrak{X}}^{\bdot} :=
\ul{\Gamma}_{X_{q}} \mcal{F}_{\mfrak{X}}^{\bdot}$.
Let us write
$X^{q} / X^{q + 1} := X_{-q} / X_{-q - 1}$
and
$\mrm{N}^{q} := \mrm{N}_{-q}$, so
$\{\mrm{N}^{q} \mcal{F}_{\mfrak{X}}^{\bdot}\}$
is a decreasing filtration.

\begin{thm} \label{thm5.3}
Suppose $\opn{char} k = 0$ and $X \subset \mfrak{X}$ is a smooth
formal embedding. Then in the niveau spectral sequence converging to
$\mrm{H}^{\bdot}(X, \mcal{F}^{\bdot}_{\mfrak{X}})$,
the $E_{1}$ term is \textup{(}in the notation of \cite{ML} Chapter
\textup{XI):}
\[ E_{1}^{p,q} =
\mrm{H}^{p+q}_{X^{p} / X^{p + 1}}(X, \mcal{F}^{\bdot}_{\mfrak{X}})
\cong
\bigoplus_{x \in X^{p} / X^{p + 1}}
\mrm{H}^{q - p} \Omega^{\bdot}_{k(x)/k} , \]
and the coboundary operator is
$(-1)^{p + 1} \sum_{(x,y)} \opn{Res}_{(x,y)}$.
\end{thm}

\begin{proof}
We shall substitute indices
$(p, q) \mapsto (-q, -p)$; this puts us in the first quadrant.
Because $\mcal{F}_{\mfrak{X}}^{\bdot}$ is a complex of flasque
sheaves, one has
\[ E_{1}^{-q, -p} = \mrm{H}^{-p - q}_{X^{-q} / X^{-q + 1}}(X,
\mcal{F}_{\mfrak{X}}^{\bdot})
\cong
\bigoplus_{x \in X_{q} / X_{q-1}}
\mrm{H}^{-p}  \mcal{H}om^{\bdot}_{\mcal{O}_{\mfrak{X}}} \left(
\widehat{\Omega}^{\bdot}_{\mfrak{X}/k},
\mcal{K}_{\mfrak{X}}(x) \right) \]
(the operator $\delta$ is trivial on the $q$-skeleton).
Fix a point $x$ of dimension $q$ and let
$B := \mcal{O}_{\mfrak{X}, (x)}$. Then
$\widehat{\Omega}^{\bdot}_{\mfrak{X}/k, (x)} \cong
\Omega^{\bdot, \mrm{sep}}_{B / k}$
and by definition
\[ \mcal{H}om^{\bdot}_{\mcal{O}_{\mfrak{X}}} \left(
\widehat{\Omega}^{\bdot}_{\mfrak{X}/k},
\mcal{K}_{\mfrak{X}}(x) \right)
\cong \opn{Dual}_{B} \Omega^{\bdot, \mrm{sep}}_{B / k} . \]
Choose a coefficient field $\sigma: K = k(x) \ar  B$.
By \cite{Ye2} Theorem 8.6 there is an isomorphism of complexes
\[ \Psi_{\sigma} :
\opn{Dual}_{B} \Omega^{\bdot, \mrm{sep}}_{B / k} \iso
\opn{Dual}_{\sigma} \Omega^{\bdot, \mrm{sep}}_{B/k} =
\opn{Hom}^{\mrm{cont}}_{K; \sigma}
(\Omega^{\bdot, \opn{sep}}_{B/k}, \omega(K)) . \]
Here $\omega(K) = \Omega^{q}_{K/k}$ and
the operator on the right is $\opn{Dual}_{\sigma}(\mrm{d})$ of
\cite{Ye2} Definition 8.1.

According to \cite{Ye3} \S 3, $k \ar B$ is formally smooth; so $B$ is
a regular local ring, and hence
$B \cong K[[\, \ul{t}\, ]]$.
Put a grading on $\Omega^{\bdot}_{K [\, \ul{t}\, ] / k}$ by
declaring
$\opn{deg} t_{i} = \opn{deg} \mrm{d} t_{i} = 1$,
and let
$V_{l} \subset \Omega^{\bdot}_{K [\, \ul{t}\, ] / k}$
be the homogeneous component of degree $l$.
In particular
$V_{0} = \Omega^{\bdot}_{K / k}$.
Since $\mrm{d}$ preserves each $V_{l}$, from the definition of
$\opn{Dual}_{\sigma}(\mrm{d})$ we see that
\[ \opn{Dual}_{\sigma} \Omega^{\bdot, \mrm{sep}}_{B/k} =
\bigoplus_{l=0}^{\infty} \opn{Hom}_{K}(V_{l}, \omega(K)) \]
as complexes. Because the $K$-linear homotopy operator in the proof
of Lemma \ref{lem5.1} also preserves $V_{l}$ we get
$\mrm{H} \opn{Hom}_{K}(V_{l}, \omega(K)) = 0$
for $l \neq 0$, and hence
\begin{equation} \label {eqn5.2}
\mrm{H}^{-p} \opn{Dual}_{\sigma} \Omega^{\bdot, \mrm{sep}}_{B/k}
\cong \mrm{H}^{-p} \opn{Hom}_{K} (\Omega^{\bdot}_{K/k}, \omega(K))
\cong \mrm{H}^{q-p} \Omega^{\bdot}_{K/k}
\end{equation}
(cf.\ proof of Lemma \ref{lem4.3}).

It remains to check that the coboundary maps match up.
Given an immediate specialization $(x,y)$, choose a pair of compatible
coefficient fields
$\sigma: k(x) \ar \mcal{O}_{\mfrak{X}, (x)} = B$ and
$\tau: k(y) \ar \mcal{O}_{\mfrak{X}, (y)} = A$ (cf.\ \cite{Ye1}
Definition 4.1.5). Set
$\widehat{B} := \mcal{O}_{\mfrak{X}, (x,y)}$, so
$f : A \ar \widehat{B}$ is a morphism of BCAs. A cohomology class
$[\phi] \in
\mrm{H}^{-p} \opn{Dual}_{B} \Omega^{\bdot, \mrm{sep}}_{B / k}$
is sent under the isomorphism (\ref{eqn5.2})
to the class $[\beta]$ of some form
$\beta \in \Omega^{q-p}_{k(x)/k}$, such that $\mrm{d} \beta = 0$ and on
$\sigma(\Omega^{\bdot}_{k(x)/k}) \subset \Omega^{\bdot,
\mrm{sep}}_{B/k}$,
$\phi$ acts like left multiplication by $\beta$. So for
$\alpha \in \Omega^{p}_{k(y)/k}$,
\begin{eqnarray*}
\opn{Tr}_{A / k(y)} \opn{Tr}_{\widehat{B} / A} \phi f \tau (\alpha)
& = & \opn{Res}_{k((x,y)) / k(y); \tau}(\beta \wedge \tau(\alpha)) \\
& = & \opn{Res}_{k((x,y)) / k(y); \tau}(\beta) \wedge \alpha .
\end{eqnarray*}
This says that under the isomorphism
\[ \mrm{H}^{-p} \opn{Dual}_{A} \Omega^{\bdot, \mrm{sep}}_{A / k}
\cong \mrm{H}^{q-p-1} \Omega^{\bdot}_{k(y)/k} , \]
the class $\delta_{(x,y)}([\phi])$ is sent to
$\opn{Res}_{(x,y)}([\beta])$.
\end{proof}

\begin{rem} Theorem \ref{thm5.3}, but with
$\opn{R} \ul{\Gamma}_{X} \Omega^{\bdot}_{Y / k}$ instead of
$\mcal{F}_{\mfrak{X}}^{\bdot}$ (cf.\ Proposition \ref{prop5.1}),
was discovered by Grothendieck
(see \cite{Gr} Footnotes 8,9), and proved by Bloch-Ogus \cite{BlO}.
Our proof is completely different from that in \cite{BlO}, and in
particular we obtain the formula for the coboundary operator as a
sum of Parshin residues. On the other hand the proof in
\cite{BlO} is valid
for a general homology theory (including $l$-adic homology).
Bloch-Ogus went on to prove additional important results, such as the
degeneration of the sheafified spectral sequence
$\mcal{E}^{p,q}_{r}$ at $r=2$, for $X$ smooth.
\end{rem}

The next result is a generalization of \cite{Ha} Theorem II.3.2.
Suppose $X \subset \mfrak{Y}$ is another smooth formal embedding. By a
morphism of embeddings $f : \mfrak{X} \ar \mfrak{Y}$ we mean a
morphism of formal schemes inducing the identity on $X$.
Since $f$ is formally finite, according to Proposition \ref{prop5.5},
$\opn{Tr}_{f} : \mcal{F}_{\mfrak{X}}^{\bdot} \ar
\mcal{F}_{\mfrak{Y}}^{\bdot}$
is a map of complexes in $\mathsf{Ab}(X)$.

\begin{cor} \label{cor5.5}
Let $f : \mfrak{X} \ar \mfrak{Y}$ be a morphism of embeddings of $X$.
Then
$\opn{Tr}_{f} : \Gamma(X, \mcal{F}_{\mfrak{X}}^{\bdot}) \ar
\Gamma(X, \mcal{F}_{\mfrak{Y}}^{\bdot})$
is a quasi-isomorphism. If $g : \mfrak{X} \ar \mfrak{Y}$
is another such morphism, then
$\mrm{H}(\opn{Tr}_{f}) = \mrm{H}(\opn{Tr}_{g})$.
\end{cor}

\begin{proof}
$\opn{Tr}_{f}$ induces a map of niveau spectral sequences
$E^{p,q}_{r}(\mfrak{X}) \ar E^{p,q}_{r}(\mfrak{Y})$.
The theorem and its proof imply that these spectral
sequences coincide for $r \geq 1$, hence
$\mrm{H}^{\bdot}(\opn{Tr}_{f})$
is an isomorphism.
The other statement is proved like in \cite{Ye3} Theorem 2.7 (cf.\ next
corollary).
\end{proof}

\begin{cor} \label{cor5.1}
The $k$-module
$\mrm{H}^{q}(X, \mcal{F}_{\mfrak{X}}^{\bdot})$
is independent of the smooth formal embedding $X \subset \mfrak{X}$.
\end{cor}

\begin{proof}
As shown in \cite{Ye3}, given any two embeddings
$X \subset \mfrak{X}$ and $X \subset \mfrak{Y}$, the completion
of their product along the diagonal
$(\mfrak{X} \times_{k} \mfrak{Y})_{/X}$ is also
a smooth formal embedding of $X$, and it projects to both $\mfrak{X}$
and $\mfrak{Y}$. Therefore by Corollary \ref{cor5.5},
$\mrm{H}^{q}(X, \mcal{F}_{\mfrak{X}}^{\bdot})$ and
$\mrm{H}^{q}(X, \mcal{F}_{\mfrak{Y}}^{\bdot})$
are isomorphic. Using triple products we see this isomorphism is
canonical.
\end{proof}

\begin{rem} \label{rem5.2}
We can use Corollary \ref{cor5.1} to {\em define}
$\mrm{H}^{\mrm{DR}}_{\bdot}(X)$ if some smooth formal embedding exists.
For a definition of $\mrm{H}^{\mrm{DR}}_{\bdot}(X)$
in general, using a system of local embeddings, see \cite{Ye3}
(cf.\ \cite{Ha} pp.\ 28-29).
\end{rem}

\begin{rem} \label{rem5.3}
In \cite{Ye4} it is shown that $\mcal{F}^{\bdot}_{\mfrak{X}}$
is naturally a DG module over the adele-De Rham complex
$\mcal{A}^{\bdot}_{\mfrak{X}} =
\ul{\mathbb{A}}^{\bdot}_{\mrm{red}}(\widehat{\Omega}^{\bdot}_{X/k})$,
and this multiplication induces
the cap product of
$\mrm{H}_{\mrm{DR}}^{\bdot}(X)$ on
$\mrm{H}^{\mrm{DR}}_{\bdot}(X)$.
\end{rem}

The next result is new (cf.\ \cite{BlO} Example 2.2):

\begin{thm} \label{thm5.1}
De Rham homology $\mrm{H}^{\mrm{DR}}_{\bdot}(-)$
is contravariant w.r.t.\ \'{e}tale morphisms.
\end{thm}

\begin{proof}
The ``topological invariance of \'{e}tale morphisms'' (see
\cite{Mi} Theorem I.3.23) implies that
the smooth formal embedding $X \subset \mfrak{X}$
induces an ``embedding of \'{e}tale sites''
$X_{\mrm{et}} \subset \mfrak{X}_{\mrm{et}}$.
By this we mean that for every \'{e}tale morphism $U \ar X$ there
is some  \'{e}tale morphism $\mfrak{U} \ar \mfrak{X}$, unique up to
isomorphism, s.t.\ $U \cong \mfrak{U} \times_{\mfrak{X}} X$
(see \cite{Ye3}). Then $U \subset \mfrak{U}$ is a smooth formal
embedding. From Proposition \ref{prop5.5} we see there is a complex
of sheaves
$\mcal{F}^{\bdot}_{\mfrak{X}_{\mrm{et}}}$ on $X_{\mrm{et}}$
with
$\mcal{F}^{\bdot}_{\mfrak{X}_{\mrm{et}}} |_{U} \cong
\mcal{F}^{\bdot}_{\mfrak{U}} \cong g^{*}
\mcal{F}^{\bdot}_{\mfrak{X}}$
for every $g : \mfrak{U} \ar \mfrak{X}$ \'{e}tale
(cf.\ \cite{Mi} Corollary II.1.6). But by Corollary \ref{cor5.1},
$\mrm{H}^{\mrm{DR}}_{\bdot}(U) =
\mrm{H}^{\bdot}(U, \mcal{F}^{\bdot}_{\mfrak{U}})$.
\end{proof}

Say $X$ is smooth irreducible of dimension $n$. Define the sheaf
$\mcal{H}^{p}_{\mrm{DR}}$
on $X_{\mrm{Zar}}$ to be the sheafification of the presheaf
$U \mapsto \mrm{H}^{p}_{\mrm{DR}}(U)$. For any point $x \in X$
let $i_{x} : \{x\} \ar X$ be the inclusion.
Let $x_{0}$ be the generic point, so
$X_{n} / X_{n - 1} = \{ x_{0} \}$.
According to \cite{BlO} there is an exact sequence of sheaves
\[ 0 \ar \mcal{H}^{p}_{\mrm{DR}} \ar
i_{x_{0}\, *} \mrm{H}^{p} \Omega^{\bdot}_{k(x_{0}) / k}
\ar \cdots \ar
\bigoplus_{x \in X_{q} / X_{q - 1}} i_{x*}
\mrm{H}^{p + q - n} \Omega^{\bdot}_{k(x)/k}
\ar \cdots \]
called the {\em arithmetic resolution}. Observe that this is a flasque
resolution.

\begin{cor} \label{cor5.2}
The coboundary operator in the arithmetic resolution of
$\mcal{H}^{p}_{\mrm{DR}}$ is
\[ (-1)^{q + 1} \sum_{(x,y)} \opn{Res}_{(x,y)} \]
where $\opn{Res}_{(x,y)}$ is the Parshin residue of Proposition
\textup{\ref{prop5.2}}.
\end{cor}

\begin{proof}
Take $\mfrak{X} = X$ in Theorem \ref{thm5.3}, and use
\cite{BlO} Theorem 4.2.
\end{proof}

\section{The Intersection Cohomology $\mcal{D}$-module of a Curve}

Suppose $Y$ is an $n$-dimensional smooth algebraic variety over
$\mbb{C}$ and $X$ is a subvariety of codimension $d$. Let
$\mcal{H}_{X}^{d} \mcal{O}_{Y}$
be the sheaf of $d$-th cohomology of $\mcal{O}_{Y}$ with support
in $X$. According to \cite{BrKa}, the holonomic $\mcal{D}_{Y}$-module
$\mcal{H}_{X}^{d} \mcal{O}_{Y}$ has a unique simple coherent
submodule $\mcal{L}(X,Y)$, and the De Rham complex
$\opn{DR} \mcal{L}(X,Y) =
\mcal{L}(X,Y) \otimes \Omega^{\bdot}_{Y^{\mrm{an}}}[n]$
is the middle perversity intersection cohomology sheaf
$\mcal{IC}^{\bdot}_{X^{\mrm{an}}}$.
Here $Y^{\mrm{an}}$ is the associated complex manifold.
The module $\mcal{L}(X,Y)$ was described explicitly using
complex-analytic methods by Vilonen \cite{Vi} and Barlet-Kashiwara
\cite{BaKa}.
These descriptions show that the fundamental class $\mrm{C}_{X/Y}$ lies
in $\mcal{L}(X,Y) \otimes \Omega^{d}_{Y/k}$, a fact proved earlier
by Kashiwara using the Riemann-Hilbert correspondence and the
decomposition theorem of Beilinson-Bernstein-Deligne-Gabber
(see \cite{Br}).

Now let $k$ be any field of characteristic $0$, $Y$ an $n$-dimensional
smooth variety over $k$, and $X \subset Y$ an integral curve with
arbitrary singularities. In this section we give a description of
$\mcal{L}(X,Y) \subset \mcal{H}_{X}^{n-1} \mcal{O}_{Y}$
in terms of algebraic residues.
As references on $\mcal{D}$-modules we suggest
\cite{Bj} and \cite{Bo} Chapter VI.

Denote by $w$ the generic point of $X$.
Pick any coefficient field
$\sigma : k(w) \ar \widehat{\mcal{O}}_{Y, w} = \mcal{O}_{Y, (w)}$.
As in \cite{Hu} Section 1 there is a residue map
\[ \opn{Res}^{\mrm{lc}}_{w, \sigma} :
\mrm{H}_{w}^{n-1} \Omega^{n}_{Y/k} \ar
\Omega^{1}_{k(w) / k}  \]
(``lc'' is for local cohomology)
defined as follows. Choose a regular system of parameters
$f_{1}, \ldots, f_{n-1}$ in $\mcal{O}_{Y, w}$,
so that
$\mcal{O}_{Y, (w)} \cong k(w)[[f_{1}, \ldots, f_{n-1}]]$.
Then for a generalized
fraction, with $\alpha \in \Omega^{1}_{k(X) / k}$, we have
\[ \opn{Res}^{\mrm{lc}}_{w, \sigma}
\gfrac{\sigma(\alpha) \wedge \mrm{d} f_{1} \wedge \cdots \wedge
\mrm{d} f_{n - 1}}
{f_{1}^{i_{1}} \cdots f_{n - 1}^{i_{n - 1}}} =
\begin{cases}
\alpha & \text{if } (i_{1}, \ldots, i_{n-1}) = (1, \ldots, 1) \\
0 & \text{otherwise} .
\end{cases} \]

Let $\pi : \tilde{X} \ar X$ be the normalization,
and let $\tilde{w}$ be the generic point of $\tilde{X}$.
For any closed point $\tilde{x} \in \tilde{X}$ the residue field
$k(\tilde{x})$ is \'{e}tale over $k$, so it lifts into
$\mcal{O}_{\tilde{X}, (\tilde{x})}$. Hence we get canonical
morphisms of BCAs
$k(\tilde{x}) \ar \mcal{O}_{\tilde{X}, (\tilde{x})} \ar
k(\tilde{w})_{(\tilde{x})}$,
and a residue map
\[ \opn{Res}_{(\tilde{w}, \tilde{x})} :
\Omega^{1}_{k(w) / k} \ar
\Omega^{1, \mrm{sep}}_{k(\tilde{w})_{(\tilde{x})} / k} \ar
k(\tilde{w}) . \]

Define
\[ \opn{Res}^{\mrm{lc}}_{(\tilde{w}, \tilde{x})} :
\mrm{H}_{w}^{n-1} \Omega^{n}_{Y/k}
\xrightarrow{\opn{Res}^{\mrm{lc}}_{w, \sigma}}
\Omega^{1}_{k(w) / k}
\xrightarrow{\opn{Res}_{(\tilde{w}, \tilde{x})}} k(\tilde{x}) . \]
We shall see later that
$\opn{Res}^{\mrm{lc}}_{(\tilde{w}, \tilde{x})}$
is independent of $\sigma$.
Note that
$\mrm{H}_{w}^{n-1} \Omega^{n}_{Y/k} =
(\mcal{H}_{X}^{n-1} \Omega^{n}_{Y/k})_{w}$.

\begin{thm} \label{thm6.6}
Let $x \in X$ be a closed point and let
$a \in (\mcal{H}_{X}^{n-1} \mcal{O}_{Y})_{x}$.
Then
$a \in \mcal{L}(X, Y)_{x}$ iff
$\opn{Res}^{\mrm{lc}}_{(\tilde{w}, \tilde{x})}(a \alpha) = 0$
for all $\alpha \in \Omega^{n}_{Y / k, x}$
and $\tilde{x} \in \pi^{-1}(x)$.
\end{thm}

This is our algebraic counterpart of Vilonen's formula in \cite{Vi}.
The proof of the theorem appears later in this section.

Fix a closed point $x \in X$.
Write $B := \mcal{O}_{Y, (w, x)}$
and
$L := \prod_{\tilde{x} \in \pi^{-1}(x)} k(\tilde{x})$.

\begin{lem}
There is a canonical morphism of BCAs $L \ar B$, and
$B \cong L((g))[[f_{1}, \ldots, f_{n - 1}]]$
for indeterminates $g, f_{1}, \ldots, f_{n - 1}$.
\end{lem}

\begin{proof}
Because $\mcal{O}_{\tilde{X}, (\tilde{x})}$ is a regular local ring
we get
$\mcal{O}_{\tilde{X}, (\tilde{x})} \cong k(\tilde{x})[[g]]$.
It is well known (cf.\ \cite{Ye1} Theorem 3.3.2) that
$k(w)_{(x)} = k(w) \otimes \mcal{O}_{X, (x)}
\cong \prod k(\tilde{w})_{(\tilde{x})}$,
hence $k(w)_{(x)} \cong L((g))$.

Choose a coefficient field $\sigma : k(w) \ar \mcal{O}_{Y, (w)}$.
It extends to a lifting
$\sigma_{(x)} : k(w)_{(x)} \ar \mcal{O}_{Y, (w, x)} = B$
(cf.\ \cite{Ye1} Lemma 3.3.9),
and $L \ar B$ is independent of $\sigma$.
Taking a system of regular parameters
$f_{1}, \ldots, f_{n - 1} \in \mcal{O}_{Y, w}$
we obtain the desired isomorphism.
\end{proof}

The BCA $A := \mcal{O}_{Y, (x)}$ is canonically an algebra over
$K := k(x)$, so
there is a morphism of BCAs
$L \otimes_{K} A \ar B$. Define a homomorphism
\[ T_{x} : \mcal{K}(B) \xrightarrow{\mrm{Tr}}
\mcal{K}(L \otimes_{K} A) \cong
L \otimes_{K} \mcal{K}(A) . \]
Since
$A \ar L \otimes_{K} A \ar B$ are topologically \'{e}tale
(relative to $k$), it follows that $T_{x}$
is a homomorphism of $\mcal{D}(A)$-modules.

Define
\[ V(x) := \opn{Coker} \left( K \ar L \right) . \]
Observe that $V(x) = 0$ iff $x$ is either a smooth point or a
geometrically unibranch singularity.
We have
$V(x)^{*} \subset L^{*}$,
where
$(-)^{*} := \opn{Hom}_{k}(-, k)$.
The isomorphism $L^{*} \cong L$ induced by $\opn{Tr}_{L / k}$
identifies
$V(x)^{*} \cong \opn{Ker}(L \xrightarrow{\opn{Tr}} K)$.

Since
$\Omega^{n}_{Y / k}[n] \ar \mcal{K}^{\bdot}_{Y}$
is a quasi-isomorphism we get a short exact sequence
\begin{equation} \label{eqn6.1}
0 \ar (\mcal{H}_{X}^{n-1} \Omega^{n}_{Y/k})_{x} \ar \mcal{K}_{Y}(w)
\xrightarrow{\delta\ } \mcal{K}_{Y}(x) \ar 0 .
\end{equation}
Also we see that
$\mcal{K}(A) = \mcal{K}_{Y}(x) \cong
\mrm{H}^{n}_{x} \Omega^{n}_{Y / k}$.
Now
$\mcal{K}_{Y}(w) = \mcal{K}(\mcal{O}_{Y, (w)}) \subset
\mcal{K}(B)$.
Because the composed map
\[ \mcal{K}_{Y}(w) \xrightarrow{T_{x}}
L \otimes_{K} \mcal{K}_{Y}(x) \xrightarrow{\opn{Tr}_{L / K} \otimes 1}
\mcal{K}_{Y}(x) \]
coincides with $\delta$,
and by the sequence (\ref{eqn6.1}), we obtain a homomorphism
of $\mcal{D}_{Y, x}$-modules
\begin{equation} \label{eqn6.11}
T_{x} : (\mcal{H}_{X}^{n-1} \Omega^{n}_{Y/k})_{x} \ar
V(x)^{*} \otimes_{K} \mrm{H}^{n}_{x} \Omega^{n}_{Y / k} .
\end{equation}

\begin{thm} \label{thm6.3}
The homomorphism $T_{x}$
induces a bijection between the lattice of nonzero
$\mcal{D}_{Y, x}$-submodules of
$(\mcal{H}_{X}^{n-1} \Omega^{n}_{Y/k})_{x}$
and the lattice of $k(x)$-submodules of $V(x)^{*}$.
\end{thm}

The proof of the theorem is given later in this section.

In order to globalize we introduce the following notation.
Let $Z$ be the reduced subscheme supported on the singular locus
$X_{\mrm{sing}}$, so $\mcal{O}_{Z} =
\prod_{x \in X_{\mrm{sing}}} k(x)$.
Then $\mcal{V} := \bigoplus_{x \in X_{\mrm{sing}}} V(x)$ and
$\mcal{H}_{Z}^{n} \mcal{O}_{Y} =
\bigoplus_{x \in X_{\mrm{sing}}} \mcal{H}_{ \{x\} }^{n} \mcal{O}_{Y}$
are $\mcal{O}_{Z}$-modules.
Using $\Omega^{n}_{Y / k} \otimes$ to switch between left and right
$\mcal{D}_{Y}$-modules, and identifying $V(x)^{*} \cong V(x)$
by the trace pairing, we see that Theorem \ref{thm6.3}
implies

\begin{cor} \label{cor6.10}
The homomorphism of $\mcal{D}_{Y}$-modules
\[ T := \sum_{x} T_{x} :
\mcal{H}_{X}^{n-1} \mcal{O}_{Y} \ar
(\mcal{H}_{Z}^{n} \mcal{O}_{Y}) \otimes_{\mcal{O}_{Z}} \mcal{V} \]
induces a bijection between the lattice of
nonzero coherent $\mcal{D}_{Y}$-submodules
of $\mcal{H}_{X}^{n-1} \mcal{O}_{Y}$
and the lattice of
$\mcal{O}_{Z}$-submodules of $\mcal{V}$.
\end{cor}

Since $\mcal{H}_{ \{x\} }^{n} \mcal{O}_{Y}$
is a simple $\mcal{D}_{Y}$-submodule, as immediate corollaries we get:

\begin{cor} \label{cor6.1}
$\mcal{H}_{X}^{n-1} \mcal{O}_{Y}$ has a unique simple coherent
$\mcal{D}_{Y}$-submodule
\blnk{2mm} \linebreak
$\mcal{L}(X,Y)$, and the sequence
\begin{equation} \label{eqn6.7}
0 \ar \mcal{L}(X,Y) \ar \mcal{H}_{X}^{n-1} \mcal{O}_{Y}
\xrightarrow{T}
(\mcal{H}_{Z}^{n} \mcal{O}_{Y}) \otimes_{\mcal{O}_{Z}} \mcal{V}  \ar 0
\end{equation}
is exact.
\end{cor}

\begin{cor} \label{cor6.2}
$\mcal{H}_{X}^{n-1} \mcal{O}_{Y}$ is a simple coherent
$\mcal{D}_{Y}$-module iff
the singularities of $X$ are all geometrically unibranch.
\end{cor}

According to Proposition \ref{prop4.7} the fundamental class
$\mrm{C}_{X/Y}$ is a double cocycle in
$\mcal{H}om(\Omega^{1}_{Y/k}, \mcal{K}^{-1}_{Y})$,
so it determines a class in
$(\mcal{H}_{X}^{n-1} \mcal{O}_{Y}) \otimes_{\mcal{O}_{Y}}
\Omega^{n-1}_{Y/k}$.

\begin{thm} \label{thm6.2}
$\mrm{C}_{X/Y} \in \mcal{L}(X,Y) \otimes_{\mcal{O}_{Y}}
\Omega^{n-1}_{Y/k}$.
\end{thm}

This of course implies that if $\alpha_{1}, \ldots, \alpha_{n}$
is a local basis of $\Omega^{n-1}_{Y/k}$ and
$\mrm{C}_{X/Y} = \sum a_{i} \otimes \alpha_{i}$, then any nonzero
$a_{i}$ generates $\mcal{L}(X,Y)$ as a $\mcal{D}_{Y}$-module.
The proof of the theorem is given later in this section.

\begin{rem}
As the referee points out, when $k = \mbb{C}$, Corollary \ref{cor6.10}
follows easily from the Riemann-Hilbert correspondence.
In that case we may consider the sheaf $\mcal{V}$ on the analytic space
$X^{\mrm{an}}$, given by
$\mcal{V} := \opn{Coker}(\mbb{C}_{X^{\mrm{an}}} \ar
\pi^{\mrm{an}}_{*} \mbb{C}_{\tilde{X}^{\mrm{an}}})$.
Now
$\mcal{IC}_{X^{\mrm{an}}} \cong \pi^{\mrm{an}}_{*}
\mbb{C}_{\tilde{X}^{\mrm{an}}}[1]$.
The triangle
$\mcal{V} \ar \mbb{C}_{X^{\mrm{an}}}[1] \ar$ \linebreak
$\mcal{IC}_{X^{\mrm{an}}} \xrightarrow{+1}$
is an exact sequence in the category of perverse sheaves, and
it is the image of (\ref{eqn6.7})
under the functor
$\opn{Sol} = \mrm{R} \mcal{H}om_{\mcal{D}_{Y^{\mrm{an}}}}((-)^{\mrm{an}},
\mcal{O}_{Y^{\mrm{an}}}[n])$.
Nonetheless ours seems to be the first purely algebraic proof
Theorem \ref{thm6.3}
and its corollaries (but cf.\ next remark).
\end{rem}

\begin{rem}
When $Y = \mbf{A}^{2}$ (i.e.\ $X$ is an affine plane curve)
and $k$ is algebraically closed,
Corollary \ref{cor6.2} was partially proved
by S.P.\ Smith \cite{Sm}, using the ring structure of $\mcal{D}(X)$.
Specifically, he proved that if $X$ has unibranch singularities, then
$\mcal{H}_{X}^{1} \mcal{O}_{Y}$ is simple.
\end{rem}

\begin{exa} \label{exa6.1}
Let $X$ be the nodal curve in
$Y = \mbf{A}^{2} = \opn{Spec} k [s, t]$
defined by
$f = s^{2} (s+1) - t^{2}$, and let $x$ be the origin.
Take
$r := t / s \in \mcal{O}_{Y, w}$, so
$s = (r+1)(r-1)$.
We see that $\tilde{X} = \opn{Spec} k[r]$ and
$r + 1, r - 1$ are regular parameters
at $\tilde{x}_{1}, \tilde{x}_{2}$ respectively on $\tilde{X}$.
For any coefficient field $\sigma$,
\[ \opn{Res}^{\mrm{lc}}_{w, \sigma}
\gfrac{\mrm{d} s \wedge \mrm{d} t}{f} =
\opn{Res}^{\mrm{lc}}_{w, \sigma}
\gfrac{- \mrm{d} (r + 1) \wedge \mrm{d} f}
{(r + 1)(r - 1) f} = \frac{- \mrm{d} (r + 1)}{(r + 1)(r - 1)}  \]
and hence
\[ \opn{Res}^{\mrm{lc}}_{(\tilde{w}, \tilde{x}_{1})}
\gfrac{\mrm{d} s \wedge \mrm{d} t}{f} =
\opn{Res}_{(\tilde{w}, \tilde{x}_{1})}
\frac{- \mrm{d} (r + 1)}{(r + 1)(r - 1)} = 2 . \]
Likewise
$\opn{Res}^{\mrm{lc}}_{(\tilde{w}, \tilde{x}_{2})}
\gfrac{\mrm{d} s \wedge \mrm{d} t}{f} = - 2$.
Therefore
$\gfrac{\mrm{d} s \wedge \mrm{d} t}{f} \notin
\mcal{L}(X, Y) \otimes \Omega^{2}_{Y / k}$.
The fundamental class is
$\mrm{C}_{X / Y} = \gfrac{\mrm{d} f}{f}$, and as generator of
$\mcal{L}(X, Y) \otimes \Omega^{2}_{Y / k}$
we may take
$\gfrac{\mrm{d} s \wedge \mrm{d} f}{f}$.
\end{exa}

Before getting to the proofs we need some general results.
Let $A$ be a BCA over $k$.
The fine topology on an $A$-module $M$ is the quotient topology
w.r.t.\ any surjection $\bigoplus A \surj M$.
The fine topology on $M$ is $k$-linear, making it a topological
$k$-module (but only a semi-topological (ST) $A$-module).
According to \cite{Ye2} Proposition 2.11.c, $A$ is a Zariski ST ring
(cf.\ ibid.\ Definition 1.7). This means that any finitely
generated $A$-module with the fine topology is separated, and any
homomorphism $M \ar N$ between such modules is topologically
strict. Furthermore if $M$ is finitely generated then it is
complete, so it is a complete linearly topologized $k$-vector space in
the sense of \cite{Ko}.

\begin{lem} \label{lem6.3}
Let $A$ be a BCA.
Suppose $M$ is a countably generated ST $A$-module with the fine
topology. Then $M$ is separated, and any submodule $M' \subset M$
is closed.
\end{lem}

\begin{proof}
Write $M = \bigcup_{i=1}^{\infty} M_{i}$ with $M_{i}$ finitely generated.
Suppose we put the fine topology on $M_{i}$. Then each
$M_{i}$ is separated and $M_{i} \ar M_{i+1}$ is strict.
By \cite{Ye1} Corollary 1.2.6 we have $M \cong \lim_{i \ar} M_{i}$
topologically, so by ibid.\ Proposition 1.1.7, $M$ is separated.
By the same token $M / M'$ is separated too, so $M'$ is closed.
\end{proof}

\begin{prop} \label{prop6.1}
Let $A \ar B$ be a morphism of BCAs, $N$ a finitely generated
$B$-module with the fine topology, and $M \subset N$ a finitely
generated $A$-module. Then
the topology on $M$ induced by $N$ equals the fine $A$-module
topology, and $M$ is closed in $N$.
\end{prop}

\begin{proof}
Since $A$ is a Zariski ST ring we may replace $M$ by any
finitely generated $A$ module $M'$, $M \subset M' \subset N$.
Therefore we can assume $N = B M$ and
$M = \bigoplus_{\mfrak{n} \in \opn{Max} B} M \cap N_{\mfrak{n}}$.
So in fact we may assume $A,B$ are both local. Like in the proof of
\cite{Ye2} Theorem 7.4 we may further assume that
$\opn{res.dim}(A \ar B) \leq 1$.

Put on $M$ the fine $A$-module topology.
Let $\bar{N}_{i} := N / \mfrak{n}^{i} N$ and
$\bar{M}_{i} := M / (M \cap \mfrak{n}^{i} N)$ with the quotient
topologies.
We claim $\bar{M}_{i} \ar \bar{N}_{i}$ is a strict monomorphism.
This is so because as $A$-modules both have the fine topology,
$\bar{M}_{i}$ is finitely generated and $\bar{N}_{i}$ is countably
generated (cf.\ part 1 in the proof of \cite{Ye1} Theorem 3.2.14).
Just as in part 2 of loc.\ cit.\ we get topological isomorphisms
$M \cong \lim_{\leftarrow i} \bar{M}_{i}$ and
$N \cong \lim_{\leftarrow i} \bar{N}_{i}$, so $M \ar N$ is a strict
monomorphism. But $M$ is complete and $N$ is separated, so $M$ must be
closed.
\end{proof}

Given a topological $k$-module $M$ we set
$M^{*} := \mrm{Hom}_{k}^{\mrm{cont}}(M,k)$
(without a topology).

\begin{lem} \label{lem6.2}
Suppose $M$ is a separated topological $k$-module. Then:
\begin{enumerate}
\item For any subset $S \subset M^{*}$ its perpendicular
$S^{\perp} \subset M$ is a closed submodule.
\item Given a closed submodule $M_{1} \subset M$, one has
$M_{1}^{\perp \perp} = M_{1}$.
\item Suppose $M_{1} \subset M_{2} \subset M$ are closed submodules.
Then there is an exact sequence \textup{(}of untopologized
$k$-modules\textup{)}
\[ 0 \ar M_{2}^{\perp} \ar M_{1}^{\perp} \ar (M_{2} / M_{1})^{*} \ar 0
. \]
\end{enumerate}
\end{lem}

\begin{proof}
See \cite{Ko} Section 10.4, 10.8.
\end{proof}

Let $M,N$ be complete separated topological $k$-modules, and
$\langle -,- \rangle : M \times N \ar k$ a continuous pairing.
We say $\langle -,- \rangle$ is a {\em topological perfect pairing}
if it induces isomorphisms $N \cong M^{*}$ and $M \cong N^{*}$.

\begin{prop} \label{prop6.3}
Assume $k \ar A$ is a morphism of BCAs. Then the residue pairing
$\langle -,- \rangle_{A/k} : A \times \mcal{K}(A) \ar k$,
$\langle a, \phi \rangle_{A/k} = \opn{Tr}_{A / k}(a \phi)$,
is a topological perfect pairing.
\end{prop}

\begin{proof}
We may assume $A$ is local. Then
$A = \lim_{\leftarrow i} A / \mfrak{m}^{i}$ and
$\mcal{K}(A) = \lim_{i \ar} \mcal{K}(A / \mfrak{m}^{i})$
topologically.
Let $K \ar A$ be a coefficient field, so both $A / \mfrak{m}^{i}$ and
$\mcal{K}(A / \mfrak{m}^{i}) \cong
\mrm{Hom}_{K}(A / \mfrak{m}^{i}, \omega(K))$
are finite $K$-modules with the fine topology. By \cite{Ye1}
Theorem 2.4.22 the pairing is perfect.
\end{proof}

From here to the end of this section we consider an integral curve $X$
embedded as a closed subscheme in a smooth irreducible $n$-dimensional
variety $Y$.
Fix a closed point $x \in X$, and set
$A := \mcal{O}_{Y, (x)}$ and $K := k(x)$.
Choosing a regular system of parameters at $x$, say
$\ul{t} = (t_{1}, \ldots, t_{n})$, allows us to write
$A = K [[\, \ul{t}\, ]]$. Let
$\mcal{D}(A) := \mrm{Diff}^{\mrm{cont}}_{A/k}(A,A)$.
Since both $K [\, \ul{t}\, ] \ar A$ and $\mcal{O}_{Y,x} \ar A$
are topologically \'{e}tale relative to $k$, we have
\[ \mcal{D}(A) \cong A \otimes_{K}
K [ \textstyle{\frac{\partial}{\partial t_{1}}}, \ldots,
\textstyle{\frac{\partial}{\partial t_{n}}}] \cong
A \otimes_{\mcal{O}_{Y,x}} \mcal{D}_{Y,x} \]
(cf.\  \cite{Ye2} Section 4).

Define $B := \mcal{O}_{Y, (w, x)}$. Since $A \ar B$ is topologically
\'{e}tale relative to $k$, we get a $k$-algebra homomorphism
$\mcal{D}(A) \ar \mcal{D}(B)$. In particular, $B$ and
$\mcal{K}(B)$ are $\mcal{D}(A)$-modules.
Define
$L := \prod_{\tilde{x} \in \pi^{-1}(x)} k(\tilde{x})$
as before.

\begin{lem} \label{lem6.1}
The multiplication map
$A \otimes_{K} L \ar B$ is injective. Its image is a
$\mcal{D}(A)$-submodule of $B$. Any $\mcal{D}(A)$-submodule of
$B$ which is finitely generated over $A$ equals $A \otimes_{K} W$
for some $K$-submodule $W \subset L$.
\end{lem}

\begin{proof}
By \cite{Kz} Proposition 8.9, if $M$ is any $\mcal{D}(A)$-module which is
finitely generated over $A$, then $M = A \otimes_{K} W$, where
$W \subset M$ is the $K$-submodule consisting of all elements killed by
the derivations $\frac{\partial}{\partial t_{i}}$.
Note that $\Omega^{1, \mrm{sep}}_{B / k}$ is free with basis
$\mrm{d} t_{1}, \ldots, \mrm{d} t_{n}$.
Thus it suffices to prove that
\begin{equation} \label{eqn6.5}
L = \{ b \in B\ |\ \ \textstyle{\frac{\partial}{\partial t_{1}}} b =
\cdots = \textstyle{\frac{\partial}{\partial t_{n}}} b = 0 \}
= \mrm{H}^{0} \Omega^{\bdot, \mrm{sep}}_{B  / k} .
\end{equation}
We know that
$B \cong L((g))[[f_{1}, \ldots, f_{n - 1}]]$,
so $B$ is topologically \'{e}tale over the polynomial algebra
$k [g, f_{1}, \ldots, f_{n - 1}]$ (relative to $k$),
and hence
$\mrm{d} g,  \mrm{d} f_{1}, \ldots$, \linebreak
$\mrm{d} f_{n - 1}$
is also a basis of $\Omega^{1, \mrm{sep}}_{B / k}$.
It follows that
$\mrm{H}^{0} \Omega^{\bdot, \mrm{sep}}_{B  / k} = L$.
\end{proof}

\begin{proof} (of Theorem \ref{thm6.3})\
Set $\mcal{M} := \mcal{H}_{X}^{n-1} \Omega^{n}_{Y/k}$
and define
$M := A \otimes_{\mcal{O}_{Y, x}} \mcal{M}_{x}$.
Tensoring the exact sequence (\ref{eqn6.1}) with $A$ we get an exact
sequence of $\mcal{D}(A)$-modules
\begin{equation} \label{eqn6.2}
0 \ar M \ar \mcal{K}(B) \xrightarrow{\delta\ } \mcal{K}(A) \ar 0 .
\end{equation}
The proof will use repeatedly the residue pairing
$\langle -,- \rangle_{B/k} : B \times \mcal{K}(B) \ar k$.
By definition of $\delta$ (cf.\ Definition \ref{dfn1.1} and \cite{Ye2}
Section 7) we see that $M = A^{\perp}$.
Consider the closed $k$-submodules
$A \subset A \otimes_{K} L \subset B$ (cf.\ Proposition \ref{prop6.1}).
Applying Lemma \ref{lem6.2} to them, and using
$V(x) = L / K$ and $\mcal{K}(A) \cong A^{*}$,
we get an exact sequence of $\mcal{D}(A)$-modules
\[ 0 \ar (A \otimes_{K} L)^{\perp} \ar M \xrightarrow{T'}
V(x)^{*} \otimes_{K} \mcal{K}(A) \ar 0
. \]
Keeping track of the operations we see that in fact
$T' = T_{x}|_{M}$.

Put the fine $A$-module topology on $M$ and $\mcal{K}(A)$, so
$M \ar V(x)^{*} \otimes_{K} \mcal{K}(A)$ is continuous.
By \cite{Ye1} Proposition 1.1.8, $\mcal{M}_{x} \ar M$ is dense.
Since $\mcal{K}(A)$ is discrete we conclude that
$\mcal{M}_{x} \ar V(x)^{*} \otimes_{K} \mcal{K}(A)$
is a surjection of $\mcal{D}_{Y, x}$-modules.
Thus any $K$-module $W \subset V(x)^{*}$ determines a distinct nonzero
$\mcal{D}_{Y, x}$-module $\mcal{N}_{x} \subset \mcal{M}_{x}$.

Conversely, say
$\mcal{N}_{x} \subset \mcal{M}_{x}$
is a nonzero $\mcal{D}_{Y, x}$-module. On any open set
$U \subset Y$ s.t.\ $U \cap X$ is smooth the module
$\mcal{M}|_{U}$ is a simple coherent $\mcal{D}_{U}$-module (by
Kashiwara's Theorem it corresponds to the $\mcal{D}_{X \cap U}$-module
$\Omega^{1}_{(X \cap U) / k}$).
Therefore the finitely generated $\mcal{D}_{Y, x}$-module $C$ defined
by
\begin{equation} \label{eqn6.4}
0 \ar \mcal{N}_{x} \ar \mcal{M}_{x} \ar C \ar 0
\end{equation}
is supported on $\{x\}$. It follows that $C \cong \mcal{K}(A)^{r}$
for some number $r$.
Tensoring  (\ref{eqn6.4}) with $A$ we get an exact sequence
of $\mcal{D}(A)$-modules
\[ 0 \ar N \ar M \ar C \ar 0  \]
with $N \subset M \subset \mcal{K}(B)$.
By faithful flatness of $\mcal{O}_{Y, x} \ar A$ we see that
$\mcal{N}_{x} = \mcal{M}_{x} \cap N$.

We put on $M,N$ the topology induced from $\mcal{K}(B)$, and on
$C$ the quotient topology from $M$.
Now $\mcal{K}(B)$ has the fine $A$-module topology and it is countably
generated over $A$ (cf.\ proof of Proposition \ref{prop6.1}), so by Lemma
\ref{lem6.3} both $M,N$ are closed in $\mcal{K}(B)$.
Using Lemma \ref{lem6.2} and the fact that
$M^{\perp} = A$ we obtain the exact sequence
\[ 0 \ar A \ar N^{\perp} \ar C^{*} \ar 0 , \]
with $N^{\perp} \subset B$.
We do not know what the topology on $C$ is; but it is a ST $A$-module.
Hence the identity map $\mcal{K}(A)^{r} \ar C$ is continuous, and it
induces an $A$-linear injection $C^{*} \ar A^{r}$.
Therefore $C^{*}$, and thus also $N^{\perp}$, are finitely generated over
$A$. According to Lemma \ref{lem6.1},
$N^{\perp} = A \otimes_{K} W$ for some $K$-module $W$,
$K \subset W \subset L$. But $N$ is closed, so $N = (N^{\perp})^{\perp}$.
\end{proof}

\begin{proof} (of Theorem \ref{thm6.2})\
For each $\tilde{x} \in \pi^{-1}(x)$ define a homomorphism
\[ T_{(\tilde{w}, \tilde{x})} : \mcal{K}(B) \xrightarrow{\mrm{Tr}}
\mcal{K}(L \otimes_{K} A) \cong
L \otimes_{K} \mcal{K}(A) \ar
k(\tilde{x}) \otimes_{K} \mcal{K}(A) , \]
so
$T_{x} = \sum T_{(\tilde{w}, \tilde{x})}$.
From the proof of Theorem \ref{thm6.3} we see that
the theorem amounts to the claim that
$T_{\tilde{x}} (\mrm{C}_{X/Y}(\alpha)) = 0$ for every
$\tilde{x}$ and $\alpha \in \Omega^{1}_{Y / k, x}$.
But $\mrm{C}_{X/Y}$ is the image of
$\mrm{C}_{X} \in \mcal{H}om(\Omega^{1}_{X/k},
\mcal{K}^{-1}_{X}(w))$,
so we can reduce our residue calculation to the
curve $\tilde{X}$. In fact it suffices to show that for every
$\alpha \in \Omega^{1}_{X / k, x}$
one has
$\opn{Res}_{(\tilde{w}, \tilde{x})} \alpha = 0$.
Since
$\alpha \in \Omega^{1}_{\tilde{X} / k, \tilde{x}}$
this is obvious.
\end{proof}

\begin{proof} (of Theorem \ref{thm6.6})\
According to \cite{SY} Corollary 0.2.11 (or \cite{Hu} Theorem 2.2)
one has
\[ \opn{Res}^{\mrm{lc}}_{(\tilde{w}, \tilde{x})} =
(1 \otimes \opn{Tr}_{A / K}) T_{(\tilde{w}, \tilde{x})} :
\mrm{H}^{n - 1}_{w} \Omega^{n}_{Y / k} \ar k(\tilde{x}) , \]
which shows that
$\opn{Res}^{\mrm{lc}}_{(\tilde{w}, \tilde{x})}$
is independent of $\sigma$. Now use Theorem \ref{thm6.3}.
\end{proof}

\begin{prob}
What is the generalization to $\opn{dim} X > 1$?
To be specific, assume $X$ has only an isolated singularity at $x$.
Then we know there is an exact sequence
\[ 0 \ar \mcal{L}(X,Y) \ar \mcal{H}_{X}^{d} \mcal{O}_{Y}
\xrightarrow{T}
\mcal{H}_{ \{x\} }^{n} \mcal{O}_{Y} \otimes_{k(x)} V(x) \ar 0  \]
for some $k(x)$-module $V(x)$.
What is the geometric data determining $V(x)$ and $T$?
Is it true that $T = \sum T_{\xi}$, a sum of ``residues''
along chains $\xi \in \pi^{-1}(x)$, for a suitable
resolution of singularities $\pi : \tilde{X} \ar X$?
\end{prob}


\end{document}